\documentclass[11pt,a4paper]{article}
\usepackage{amsfonts}
 \title{\textbf{Some remarks on the generalized Tanaka-Webster connection of a contact metric manifold}}
 \date{}
 \author{\begin{tabular}{cc}
  \textsc{Beniamino Cappelletti Montano}\\
  Department of Mathematics,  University of Bari \\
  Via E. Orabona, 4 \\
  I-70125 Bari (Italy) \\
  \textsf{cappelletti@dm.uniba.it}
   \end{tabular}}
\usepackage{graphicx}
\usepackage{amsmath}
\newtheorem{theorem}{Theorem}[section]

\newtheorem{corollary}[theorem]{Corollary}

\newtheorem{definition}[theorem]{Definition}
\newtheorem{example}[theorem]{Example}

\newtheorem{lemma}[theorem]{Lemma}

\newtheorem{proposition}[theorem]{Proposition}
\newtheorem{remark}[theorem]{Remark}

\newenvironment{proof}[1][Proof]{\textbf{#1.} }{\ \rule{0.5em}{0.5em}}

\begin{document}

\maketitle

\begin{abstract}
We find \ necessary and sufficient \ conditions \ for the
bi-Legendrian connection \ $\nabla$  associated \ to a bi-Legendrian
structure $\left(\cal F,\cal G\right)$ \ on a contact metric
manifold \ $(M^{2n+1},\phi,\xi,\eta,g)$ being a metric connection
with respect to the associated metric $g$ and then we give
conditions ensuring that $\nabla$ coincides with the (generalized)
Tanaka-Webster connection of $(M^{2n+1},\phi,\xi,\eta,g)$. Using
these results, we give some interpretations of the Tanaka-Webster
connection and we study the interplays between the Tanaka-Webster,
the bi-Legendrian and the Levi Civita connection in a Sasakian
manifold.
\end{abstract}
\textbf{2000 Mathematics Subject Classification.} 53B05, 53C12,
53C15.
\\
\textbf{Keywords and phrases.} Tanaka-Webster connection, Legendrian
foliations, bi-Legen-drian connection, contact metric structure,
Sasakian manifold.

\section{Introduction}

In this paper we study some properties of the (generalized)
Tanaka-Webster connection of a contact metric manifold
$(M^{2n+1},\phi,\xi,\eta,g)$. This connection has been introduced by
S. Tanno (cf. \cite{tanno}) as a generalization of the well-known
connection defined at the end of the 70's by N. Tanaka in
\cite{tanaka} and, independently, by S. M. Webster in
\cite{webster}, in the context of CR-geometry.  We put in relation
the (generalized) Tanaka-Webster connection with the theory of
Legendrian foliations on contact metric manifolds (cf.
\cite{jayne1}, \cite{libermann}, \cite{pang}). In particular, in
\cite{mino2} the author has attached to any Legendrian foliation a
canonical connection, called \emph{bi-Legendrian connection}, and in
\cite{cappellettimontano1} he has found many applications of this
connection in the theory of Legendrian foliations. In this paper we
find conditions for which the Tanaka-Webster connection and the
bi-Legendrian connection associated to a given Legendrian foliation
coincide. We discuss some consequences of these results and give new
interpretations both of Tanaka-Webster and of bi-Legendrian
connections. For the latter, more precisely, we prove that the
bi-Legendrian connection associated to a given Legendrian foliation
on a contact manifold $(M^{2n+1},\eta)$ can be viewed as the
Tanaka-Webster connection of a suitable Sasakian structure
$\left(\phi,\xi,\eta,g\right)$ on $M^{2n+1}$ and they are
\emph{contact metric connections} in the sense of \cite{nicolaescu}.
From this and other theorems which we will prove in $\S$
\ref{confronto} and $\S$ \ref{interpretazione}, compared with the
analogous results in even dimension, we see that the Tanaka-Webster
connection of a Sasakian manifold plays the role of the Levi Civita
connection on a K\"{a}hlerian manifold. Finally, in $\S$
\ref{esempi}, we present some examples and counterexamples, for
instance we construct a Sasakian structure on $S^3$, endowed with a
non-flat bi-Legendrian structure, for which the Tanaka-Webster
connection and the bi-Legendrian connection do not coincide.

The framework of this paper are contact metric manifolds. Recall
that a \emph{contact structure} on an odd dimensional smooth
manifold $M^{2n+1}$ is given by a $1$-form $\eta$ satisfying  \
$\eta\wedge\left(d\eta\right)^n\neq 0$ everywhere on $M^{2n+1}$. It
is well known that given $\eta$ there exists a unique vector field
$\xi$, called \emph{Reeb vector field}, such that
$d\eta\left(\xi,\cdot\right)=0$ and $\eta\left(\xi\right)=1$. The
distribution defined by $\ker\left(\eta\right)$ is called the
\emph{contact distribution} and is denoted by $\cal D$. Then the
tangent bundle of $M^{2n+1}$ splits as the direct sum
$TM^{2n+1}=\cal D\oplus\mathbb{R}\xi$. A Riemannian metric $g$ is an
\emph{associated metric} for a contact form $\eta$ if the following
two conditions hold:
\begin{description}
    \item[(i)] $g\left(V,\xi\right)=\eta\left(V\right)$ for all $V\in\Gamma(TM^{2n+1})$, that is
    $\xi$ is orthogonal to $\cal D$;
    \item[(ii)] there exists a tensor field $\phi$ of type $(1,1)$ on
    $M^{2n+1}$ such that $\phi^2=-I+\eta\otimes\xi$ and
    $d\eta\left(V,W\right)=g\left(V,\phi W\right)$ for all
    $V,W\in\Gamma(TM^{2n+1})$.
\end{description}
Moreover, from (i) and (ii) one can prove the following well-known
relations (cf. \cite{blair1}):
\begin{equation*}
\phi\xi=0, \ \eta\circ\phi=0, \ g(\phi V,\phi
W)=g(V,W)-\eta(V)\eta(W), \ g(\phi V,W)=-g(V,\phi W)
\end{equation*}
for all $V,W\in\Gamma(TM^{2n+1})$. We refer to
$\left(\phi,\xi,\eta,g\right)$ as a \emph{contact metric structure}
and to $M^{2n+1}$ with such a structure as a \emph{contact metric
manifold}. A contact metric manifold $(M^{2n+1},\phi,\xi,\eta,g)$ is
called a \emph{Sasakian manifold} if it is \emph{normal}, i.e. if
the tensor field $N:=\left[\phi,\phi\right]+2d\eta\otimes\xi$
vanishes identically. In terms of the covariant derivative of $\phi$
the Sasakian condition is
\begin{equation}\label{condizionesasaki}
(\hat{\nabla}_{V}\phi)W=g\left(V,W\right)\xi-\eta\left(W\right)V,
\end{equation}
where $\hat{\nabla}$ denotes, and will denote in all this paper, the
Levi Civita connection. In the study of contact metric manifolds it
is useful to define a tensor field $h$ by $h=\frac{1}{2}{\cal
L}_{\xi}\phi$. The operator $h$ is symmetric, anti-commutes with
$\phi$, satisfies $h\xi=0$ and it vanishes if and only if $\xi$ is a
Killing vector field (in this case the contact metric manifold in
question is said to be \emph{K-contact}; it is easy to show that a
Sasakian manifold is also K-contact). Moreover,
\begin{equation}\label{condizione}
\hat{\nabla}_{V}\xi=-\phi V-\phi h V
\end{equation}
holds for all $V\in\Gamma(TM^{2n+1})$. For the proofs of all these
properties and more details on contact metric manifolds we refer the
reader to \cite{blair1}.

Given a contact metric manifold
$(M^{2n+1},\phi,\xi,\eta,g)$, there is defined on
$M^{2n+1}$ a canonical connection, called the \emph{generalized
Tanaka-Webster connection} or, simply, the Tanaka-Webster
connection of the contact metric manifold
$(M^{2n+1},\phi,\xi,\eta,g)$. This connection is
defined by the following formula:
\begin{equation}\label{tanakaconnection}
    ^{\ast}\nabla_{V}W:=\hat{\nabla}_{V}W+\eta\left(V\right)\phi
    W+\eta\left(W\right)\left(\phi V+\phi h
    V\right)+d\eta\left(V+h V,W\right)\xi,
\end{equation}
for all $V,W\in\Gamma(TM^{2n+1})$. The torsion tensor of this
connection has the following expression:
\begin{equation}\label{tanakatorsion}
^{\ast}T\left(V,W\right)=\eta\left(W\right)\phi h V -
\eta\left(V\right)\phi h W + 2g\left(V,\phi W\right)\xi.
\end{equation}
Tanno (\cite{tanno}) found a characterization of this connection. He
proved that the Tanaka-Webster connection $^{\ast}\nabla$ is the
unique linear connection on $M^{2n+1}$ such that
\begin{align}\label{listatanaka}
\nonumber  \textrm{(i)} \ \ & ^\ast\nabla g=0, \
^\ast\nabla\eta=0, \
  ^\ast\nabla\xi=0,\\
  \textrm{(ii)} \ \ &
  \left(^\ast\nabla_V\phi\right)W=(\hat{\nabla}_V\phi)W-g\left(V+hV,W\right)\xi+\eta\left(W\right)\left(V+hV\right),\\
\nonumber  \textrm{(iii)} \ \ & ^\ast T\left(\xi,\phi
V\right)=-\phi{^\ast
  T}\left(\xi,V\right),\\
\nonumber  \textrm{(iv)} \ \ &
  ^{\ast}T\left(Z,Z'\right)=2d\eta\left(Z,Z'\right)\xi \
    \textrm{ for all }  Z,Z'\in\Gamma\left(\cal{D}\right).
\end{align}
This connection agrees with the connection of Tanaka in
\cite{tanaka} when the contact metric manifold is a strongly
pseudo-convex (integrable) CR-manifold.

\medskip

All manifolds considered here are assumed to be smooth i.e.  of the
class $\mathcal{C}^{\infty}$, and connected; we denote by
$\Gamma(\cdot)$ the set of all sections of a corresponding bundle.
We use the convention that $2u\wedge v=u\otimes v-v\otimes u$.

\section{Bi-Legendrian connections}\label{preliminari}

The contact condition $\eta\wedge\left(d\eta\right)^n\neq 0$ can be
interpreted geometrically saying that the contact distribution is
far from being integrable as possible. One can prove that the
maximal dimension of an integrable subbundle $L$ of $\cal D$ is $n$.
In this case $L$ necessarily satisfies the condition $d\eta(X,X')=0$
for all $X,X'\in\Gamma(L)$, since
$2d\eta(X,X')=X(\eta(X'))-X'(\eta(X))-\eta([X,X'])=0$, $L$ being
integrable. This motivates the following definition.

\begin{definition}
A \emph{Legendrian distribution} on a contact manifold
$(M^{2n+1},\eta)$ is an $n$-dimensional subbundle $L$ of the contact
distribution such that $d\eta\left(X,X'\right)=0$ for all
$X,X'\in\Gamma\left(L\right)$. When $L$ is integrable, it defines a
\emph{Legendrian foliation} of $(M^{2n+1},\eta)$. Equivalently, a
Legendrian foliation of $(M^{2n+1},\eta)$ is a foliation of
$M^{2n+1}$ whose leaves are $n$-dimensional $C$-totally real
submanifolds of $(M^{2n+1},\eta)$.
\end{definition}

Legendrian foliations have been extensively investigated in recent
years from various points of views (cf. \cite{jayne1},
\cite{libermann}, \cite{pang}). In particular M. Y. Pang provided a
classification of Legendrian foliations by means of a bilinear
symmetric form $\Pi_{\mathcal F}$ on the tangent bundle of the
foliation, defined by $\Pi_{\mathcal
F}\left(X,X'\right)=-\left({\mathcal L}_{X}{\mathcal
L}_{X'}\eta\right)\left(\xi\right)$. He called a Legendrian
foliation $\mathcal F$  \emph{non-degenerate}, \emph{degenerate} or
\emph{flat} according to the circumstance that the bilinear form
$\Pi_{\mathcal F}$ is non-degenerate, degenerate or vanishes
identically, respectively. A geometrical interpretation of this
classification is given in the following lemma.

\begin{lemma}[\cite{jayne1}]\label{classificazione}
Let $\left(M^{2n+1},\phi,\xi,\eta,g\right)$ be a contact metric
manifold foliated by a Legendrian foliation $\mathcal F$. Then
\begin{enumerate}
    \item[\emph{(a)}] $\mathcal F$ is flat if and only if
    $\left[\xi,X\right]\in\Gamma\left(T{\mathcal F}\right)$ for all
    $X\in\Gamma\left(T{\mathcal F}\right)$,
    \item[\emph{(b)}] $\mathcal F$ is degenerate if and only if there exist
    $X\in\Gamma\left(T{\mathcal F}\right)$ such that $\left[\xi,X\right]\in\Gamma\left(T{\mathcal F}\right)$,
    \item[\emph{(c)}] $\mathcal F$ is non-degenerate if and only if
    $\left[\xi,X\right]\notin \Gamma\left(T{\mathcal F}\right)$ for all
    $X\in\Gamma\left(T{\mathcal F}\right)$.
\end{enumerate}
\end{lemma}

Lemma \ref{classificazione} also allows us to extend the notion
non-degenerateness, degenerateness and flatness to Legendrian
distributions. Thus we say that a Legendrian distribution $L$ is
\emph{flat} if $\left[\xi,X\right]\in\Gamma\left(L\right)$ for all
$X\in\Gamma\left(L\right)$, \emph{degenerate} if there exist
$X\in\Gamma\left(L\right)$ such that
$\left[\xi,X\right]\in\Gamma\left(L\right)$, and
\emph{non-degenerate} if $\left[\xi,X\right]\notin
\Gamma\left(L\right)$ for all $X\in\Gamma\left(L\right)$.
\medskip

By an \emph{almost bi-Legendrian manifold} we mean a contact
manifold $(M^{2n+1},\eta)$ endowed with two transversal Legendrian
distributions $L_1$ and $L_2$. Thus, in particular, the tangent
bundle of $M^{2n+1}$ splits up as the direct sum
$TM^{2n+1}=L_1\oplus L_2\oplus\mathbb{R}\xi$. When both $L_1$ and
$L_2$ are integrable we speak of \emph{bi-Legendrian manifold}
(\cite{mino2}). An (almost) bi-Legendrian manifold is said to be
flat, degenerate or non-degenerate if and only if both the
Legendrian distributions are flat, degenerate or non-degenerate,
respectively.

In \cite{mino2} it has been attached to any almost bi-Legendrian
manifold a canonical connection which plays an important role in the
study of almost bi-Legendrian manifolds.

\begin{theorem}[\cite{mino2}]\label{biconnection}
Let $(M^{2n+1},\eta,L_1,L_2)$ be an almost bi-Legendrian manifold.
There exists a unique linear connection ${\nabla}$ on $M^{2n+1}$
such that
\begin{align}\label{lista}
 \nonumber \emph{(i)} \ \ &\nabla L_1\subset L_1, \ \nabla L_2\subset L_2, \ \nabla\left(\mathbb{R}\xi\right)\subset\mathbb{R}\xi;\\
  \emph{(ii)} \ \ &\nabla d\eta=0;\\
 \nonumber \emph{(iii)} \ \
  &T\left(X,Y\right)=2d\eta\left(X,Y\right){\xi}, \ \textrm{ for all }
  X\in\Gamma\left(L_1\right), Y\in\Gamma\left(L_2\right),\\
 \nonumber
 &T\left(V,\xi\right)=\left[\xi,V_{L_1}\right]_{L_2}+\left[\xi,V_{L_2}\right]_{L_1},
 \  \textrm{ for all } V\in\Gamma(TM^{2n+1}),
\end{align}
where ${T}$ denotes the torsion tensor of ${\nabla}$ and $X_{L_1}$
and $X_{L_2}$ the projections of $X$ onto the subbundles $L_1$ and
$L_2$ of $TM^{2n+1}$, respectively.
\end{theorem}

Such  a  connection  is  called  the  \emph{bi-Legendrian
connection}  of   the  almost bi-Legendrian manifold
$(M^{2n+1},\eta,L_1,L_2)$. We recall the explicit construction of
this connection. First, for any two vector fields $V$ and $W$ on
$M^{2n+1}$, let $H(V,W)$ be the unique section of $\mathcal D$ such
that
\begin{equation}\
i_{H\left(V,W\right)}d\eta|_{\mathcal{D}}=\left({\mathcal{L}}_{V}i_{W}d\eta\right)|_{\mathcal{D}},
\end{equation}
that is, for every $Z\in\Gamma({\mathcal{D}})$,
$d\eta\left(H\left(V,W\right),Z\right)=V\left(d\eta\left(W,Z\right)\right)-d\eta\left(W,\left[V,Z\right]\right)$.
The existence and the uniqueness of this section depends on the fact
that the 2-form $d\eta$ is non-degenerate on ${\mathcal{D}}$. The
main properties of the operator $H$ are collected in the following
lemma.

\begin{lemma}[\cite{mino2}]\label{proprietaH}
For every $f\in C^{\infty}(M^{2n+1})$ and
$V,V^{\prime},W,W^{\prime}\in\Gamma(TM^{2n+1})$ we have:
\begin{enumerate}
    \item
$H\left(V+V^{\prime},W\right)=H\left(V,W\right)+H\left(V^{\prime},W\right)$,
\
$H\left(V,W+W^{\prime}\right)=H\left(V,W\right)+H\left(V,W^{\prime}\right)$
    \item
$H\left(V,fW\right)=fH\left(V,W\right)+V\left(f\right)W_{\mathcal{D}}$
    \item $H\left(fV,W\right)=fH\left(V,W\right)$, \ if \
    $d\eta\left(V,W\right)=0$,
    \item $H(V,\xi)=0$, \ $H(\xi,W)=[\xi,W]_{\mathcal D}$,
\end{enumerate}
where $W_{\cal{D}}$ denotes the projection of $W$ onto the subbundle
$\mathcal{D}$ of $TM^{2n+1}$.
\end{lemma}

Using Lemma \ref{proprietaH} one can define a connection
$\nabla^{L_1}$ on the bundle $L_1$ and a connection $\nabla^{L_2}$
on the bundle $L_2$ setting for all $W\in
\Gamma\left(TM^{2n+1}\right)$, $X\in \Gamma(L_1)$, $Y\in\Gamma(L_2)$
\begin{gather*}
\nabla^{L_1}_{W}X:=H\left(W_{L_1},X\right)_{L_1}+\left[W_{L_2},X\right]_{L_1}+\left[W_{\mathbb{R}\xi},X\right]_{L_1},\\
\nabla^{L_2}_{W}Y:=H\left(W_{L_2},Y\right)_{L_2}+\left[W_{L_1},Y\right]_{L_2}+\left[W_{\mathbb{R}\xi},Y\right]_{L_2}.
\end{gather*}
Moreover, we define a connection $\nabla^{\mathbb{R}\xi}$ on the
line bundle $\mathbb{R}\xi$ requiring that
$\nabla^{\mathbb{R}\xi}\xi=0$, thus setting
\begin{equation*}
\nabla^{\mathbb{R}\xi}_{W}Z:=W\left(\eta\left(Z\right)\right)\xi
\end{equation*}
for all $Z\in\Gamma(\mathbb{R}\xi)$. Then from $\nabla^{L_1}$,
$\nabla^{L_2}$ and $\nabla^{\mathbb{R}\xi}$ one can define a global
connection on $M$ by putting for any $V,W \in\Gamma(TM^{2n+1})$,
\begin{equation*}
\nabla_{W}V:=\nabla^{L_1}_{W}V_{L_1}+\nabla^{L_2}_{W}V_{L_2}+\nabla^{\mathbb{R}\xi}_{W}V_{\mathbb{R}\xi}.
\end{equation*}
In particular it follows that, for all $W\in\Gamma(TM^{2n+1})$,
$\nabla_{\xi}W=[\xi,W_{L_1}]_{L_1}+[\xi,W_{L_2}]_{L_2}+\xi(\eta(W))\xi$.
The above connection is called the bi-Legendrian connection
associated to the almost bi-Legendrian manifold
$(M^{2n+1},\eta,L_1,L_2)$ and it can be characterized as the unique
linear connection on $M$ satisfying \eqref{lista}. Further
properties of this connection are collected in the following
propositions.

\begin{proposition}[\cite{mino2}]\label{torsion}
The torsion tensor field of the bi-Legendrian connection of an
almost bi-Legendrian manifold $(M^{2n+1},\eta,L_1,L_2)$ is given by
\begin{enumerate}
  \item[\emph{(i)}] $T(X,X')=-[X,X']_{L_2}$ \
for $X,X'\in\Gamma(L_1)$,
  \item[\emph{(ii)}] $T(Y,Y')=-[Y,Y']_{L_1}$ \
for $Y,Y'\in\Gamma(L_2)$,
  \item[\emph{(iii)}] $T\left(X,Y\right)=2d\eta\left(X,Y\right)\xi$ \ for
$X\in\Gamma(L_1)$, $Y\in\Gamma(L_2)$,
  \item[\emph{(iv)}]
  $T\left(W,\xi\right)=[\xi,W_{L_1}]_{L_2}+[\xi,W_{L_2}]_{L_1}$ \
for $W\in\Gamma(TM^{2n+1})$.
\end{enumerate}
In particular, if $L_1$ and $L_2$ are flat then the terms $T(W,\xi)$
vanish, and if $L_1$ and $L_2$ are integrable then $\nabla$ is
torsion free along the leaves of the Legendrian foliations defined
by $L_1$ and $L_2$.
\end{proposition}

\begin{proposition}[\cite{cappellettimontano1}]\label{proprieta}
Let $(M^{2n+1},\eta,L_1,L_2)$ be an almost bi-Legendrian manifold
and let ${\nabla}$ denote the corresponding bi-Legendrian
connection. Then the $1$-form $\eta$ is ${\nabla}$-parallel, the
parallel transport along curves preserves the distributions $L_1$
and $L_2$, and  if $L_1, L_2$ are integrable and flat the curvature
tensor field of $\nabla$ vanishes along the leaves of the foliations
defined by $L_1$, $L_1\oplus\mathbb{R}\xi$, $L_2$ and
$L_2\oplus\mathbb{R}\xi$.
\end{proposition}

Proposition \ref{proprieta} gives a further geometrical
interpretation of the flatness of a bi-Legendrian structure. It
implies that the leaves of the Legendrian foliations in question
admit a canonical flat affine structure. This always holds in
symplectic geometry: for a symplectic manifold foliated by a
Lagrangian foliation Weinstein proved that each leaf possesses a
natural flat connection; moreover, in case the symplectic manifold
in question admits two transversal Lagrangian foliations, Hess
proved that this connection extends to a symplectic connection on
the tangent bundle called bi-Lagrangian connection (see the Appendix
for more details).

Thus the flatness of a bi-Legendrian structure, and more in general
of a Legendrian distribution, seems to be a quite natural condition
for comparing  Legendrian and  Lagrangian foliations. This can be
seen also in the following results.

\begin{proposition}
Every contact manifold $(M^{2n+1},\eta)$ endowed with a Legendrian
distribution $L$ embeds into a symplectic manifold $(C,\omega)$
endowed with a Lagrangian distribution $L^{C}$. Furthermore $L^{C}$
is integrable if and only if $L$ is integrable and flat.
\end{proposition}
\begin{proof}
Let $C=M^{2n+1}\times\mathbb{R}^{+}$ be the cone on $M^{2n+1}$ and
$\omega$ be the symplectic form on $C$ defined by
$\omega=e^{t}(d\eta-\eta\wedge dt)=d\lambda$,  $\lambda=e^t\eta$. We
set $L^{C}:=L\oplus\mathbb{R}\xi$, considered as a distribution on
$C$. Then for all $X,X'\in\Gamma(L)$ we have $\omega(X,X')=e^t
d\eta(X,X')=0$ and $\omega(X,\xi)=e^t d\eta(X,\xi)=0$, from which it
follows that $L^{C}$ is Lagrangian. The final part of the statement
is then a direct consequence of the definition of $L^{C}$.
\end{proof}

\begin{theorem}[\cite{cappellettimontano1}]
Let $(M^{2n+1},\eta)$ be a regular contact manifold endowed with a
Lagrangian distribution $L$. Then $L$ projects onto a Lagrangian
distribution on the space of leaves of $M^{2n+1}$ by the
$1$-dimensional foliation defined by $\xi$ if and only if $L$ is
flat. Furthermore, if $M^{2n+1}$ is endowed with a flat almost
bi-Legendrian structure $(L_1,L_2)$, the bi-Legendrian connection
associated to $(L_1,L_2)$ projects to the bi-Lagrangian connection
associated to the projection of $(L_1,L_2)$ on the space of leaves
of $M^{2n+1}$.
\end{theorem}

On the other hand the flatness of a Legendrian foliations implies
also some strong topological obstructions, such as a vanishing
phenomenon for the characteristic classes
(\cite{cappellettimontano1}). Moreover, we remark that there are
also several examples of non-flat Legendrian foliations (see for
instance the following Example \ref{kappamu}).

\bigskip

Any contact manifold $(M^{2n+1},\eta)$ endowed with a Legendrian
distribution $L$ admits a canonical almost bi-Legendrian structure.
Indeed let $(\phi,\xi,\eta,g)$ be a compatible contact metric
structure. Then from the relation $d\eta(\phi V,\phi W)=d\eta(V,W)$
it easily follows that $Q:=\phi L$ is a Legendrian distribution on
$M^{2n+1}$ which is orthogonal to $L$. Thus the tangent bundle of
$M^{2n+1}$ splits as the orthogonal sum  $TM^{2n+1}=L\oplus
Q\oplus\mathbb{R}\xi$. $Q$ is called the \emph{conjugate Legendrian
distribution} of $L$, and in general is not integrable even if $L$
is. Some conditions ensuring the integrability of the conjugate
Legendrian distribution of a Legendrian foliation of a contact
metric manifold are given in \cite{jayne1}.

In this article we mainly study the bi-Legendrian connection
$\nabla$ associated to the almost bi-Legendrian structure $(L,Q)$,
with $Q=\phi L$, on a contact metric manifold
$(M^{2n+1},\phi,\xi,\eta,g)$. We start finding conditions ensuring
that $\nabla$ is a metric connection with respect to the associated
metric $g$.

\begin{proposition}\label{metrica}
Let $(M^{2n+1},\phi,\xi,\eta,g)$ be a contact metric manifold
endowed with a Legendrian distribution $L$. Let $Q:=\phi L$ be the
conjugate Legendrian distribution of $L$ and $\nabla$ the
bi-Legendrian connection associated to $\left(L,Q\right)$. Then the
following statements are equivalent:
\begin{description}
    \item[(i)] $\nabla g=0$;
    \item[(ii)] $\nabla\phi=0$;
    \item[(iii)] $\nabla_{X}X'=-\left(\phi\left[X,\phi
    X'\right]\right)_L$ for all $X,X'\in\Gamma\left(L\right)$, $\nabla_{Y}Y'=-\left(\phi\left[Y,\phi
    Y'\right]\right)_Q$ for all $Y,Y'\in\Gamma\left(Q\right)$ and
    the tensor $h$ maps the subbundle $L$
    onto $L$ and the subbundle $Q$ onto $Q$;
    \item[(iv)] $g$ is a bundle-like metric with respect
both to the distribution $L\oplus \mathbb{R}\xi$ and to the
distribution $Q\oplus \mathbb{R}\xi$.
\end{description}
Furthermore, assuming $L$ and $Q$ integrable, (i)--(iv) are
equivalent to the total geodesicity  of the Legendrian foliations
defined by $L$ and $Q$
\end{proposition}
\begin{proof}
In order to prove the equivalence of (i), (ii), (iii) and (iv) it is
sufficient to prove the following implications: (i) $\Rightarrow$
(ii) $\Rightarrow$ (iii) $\Rightarrow$ (iv) $\Rightarrow$ (i).\\
(i) $\Rightarrow$ (ii) \ \ Since $d\eta$ is $\nabla$-parallel and
$d\eta(\cdot,\cdot)=g(\cdot,\phi\cdot)$,  under the assumption that
the bi-Legendrian connection is metric, we have easily that
\begin{align*}
0=(\nabla_{V}d\eta)(W,W')=(\nabla_{V}g)(W,\phi
W')+g(W,(\nabla_{V}\phi)W')=g(W,(\nabla_{V}\phi)W').
\end{align*}
for all $V,W,W'\in\Gamma(TM^{2n+1})$, from which (ii) holds.\\
(ii) $\Rightarrow$ (iii) \ \ Assuming $\nabla\phi=0$, it follows
that, for all $X,X'\in\Gamma(L)$,
$0=\left(\nabla_X\phi\right)X'=\nabla_{X}\phi X'-\phi\nabla_{X}X'$,
from which, applying $\phi$ and tacking into account that
$\nabla\eta=0$, we get
$\nabla_XX'=\eta(\nabla_{X}X')\xi-\phi\nabla_{X}\phi
X'=X(\eta(X'))\xi-\phi(\left[X,\phi
X'\right]_Q)=-\left(\phi\left[X,\phi X'\right]\right)_L$. In the
same way one finds $\nabla_YY'=-\left(\phi\left[Y,\phi
Y'\right]\right)_Q$ for all $Y,Y'\in\Gamma\left(Q\right)$. Next, for
all $X\in\Gamma\left(L\right)$, we have
\begin{equation}\label{formulah}
2\left(h X\right)_Q=\left[\xi,\phi
X\right]_Q-\phi\left(\left[\xi,X\right]_L\right)=\nabla_{\xi}\phi
X-\phi\nabla_{\xi}X=\left(\nabla_{\xi}\phi\right)X=0,
\end{equation}
and, analogously,  $2\left(h
Y\right)_L=\left(\nabla_{\xi}\phi\right)Y=0$
for all $Y\in\Gamma\left(Q\right)$.\\
(iii) $\Rightarrow$ (iv) \ \ Let us suppose that (iii) holds. Then
for all $X,X',X''\in\Gamma(L)$ we have
\begin{align*}
({\mathcal L}_{X}g)(\phi X',\phi X'')&=X(g(\phi X',\phi
X''))-g([X,\phi X'],\phi X'')-g(\phi X',[X,\phi
X''])\\
&=X(g(\phi X',\phi X''))-g([X,\phi X']_Q,\phi X'')+g(X',(\phi[X,\phi
X''])_L)\\
&=X(g(\phi X',\phi
X''))-g(\nabla_{X}\phi X',\phi X'')-g(X',\nabla_{X}X'')\\
&=X(d\eta(\phi X',X''))-d\eta(\nabla_{X}\phi X',X'')-d\eta(\phi X',\nabla_{X}X'')\\
&=(\nabla_{X}d\eta)(\phi X',X'')=0,
\end{align*}
since $d\eta$ is $\nabla$-parallel. Next, note that by
\eqref{formulah} we get $(\nabla_{\xi}\phi)X=2(hX)_{Q}=0$ for all
$X\in\Gamma(L)$ and, analogously, $(\nabla_{\xi}\phi)Y=2(hY)_{L}=0$
for all $Y\in\Gamma(Q)$. Using this we have, for all
$X',X''\in\Gamma(L)$,
\begin{align*}
({\mathcal L}_{\xi}g)(\phi X',\phi X'')&=\xi(g(\phi X',\phi
X''))-g([\xi,\phi X']_Q,\phi X'')-g(\phi X',[\xi,\phi X'']_Q)\\
&=\xi(g(\phi X',\phi X''))-g(\nabla_{\xi}\phi X',\phi X'')-g(\phi
X',\nabla_{\xi}\phi X'')\\
&=\xi(g(\phi X',\phi X''))-g(\nabla_{\xi}\phi X',\phi X'')-g(\phi
X',\phi\nabla_{\xi}X'')\\
&=(\nabla_{\xi}d\eta)(\phi X',X'')=0.
\end{align*}
Arguing in a similar way one can prove that $({\mathcal
L}_{Y}g)(X',X'')=0$ and $({\mathcal L}_{\xi}g)(X',X'')=0$ for all
$Y\in\Gamma(Q)$ and $X',X''\in\Gamma(L)$.\\
(iv) $\Rightarrow$ (i) \ \ Since the bi-Legendrian connection
$\nabla$ preserves the orthogonal decomposition $TM^{2n+1}=L\oplus
Q\oplus\mathbb{R}\xi$, in order to prove that $\nabla$ is metric it
is enough to check that $(\nabla_{V}g)(X',X'')=0$ and
$(\nabla_{V}g)(Y',Y'')=0$ for all $V\in\Gamma(TM^{2n+1})$,
$X',X''\in\Gamma(L)$ and $Y',Y''\in\Gamma(Q)$. Using (iv) we get
\begin{align}\label{formulautile}
(\nabla_{X}g)(X',X'')&=X(g(X',X''))-g(\nabla_{X}X',X'')-g(X',\nabla_{X}X'')\nonumber\\
&=X(g(X',X''))-g(H(X,X')_L,X'')-g(X',H(X,X'')_{L})\nonumber\\
&=-X(d\eta(X',\phi X''))+d\eta(H(X,X'),\phi X'')+d\eta(H(X,X''),\phi X')\nonumber\\
&=-X(d\eta(X',\phi X''))+X(d\eta(X',\phi X''))-d\eta(X',[X,\phi
X''])\nonumber\\
&\quad+X(d\eta(X'',\phi X'))-d\eta(X'',[X,\phi X'])\\
&=-X(g(X',X''))-g(X',\phi[X,\phi X''])-g(X'',\phi[X,\phi X'])\nonumber\\
&=-X(g(\phi X',\phi X''))+g([X,\phi X'],\phi X'')+g(\phi X',[X,\phi
X''])\nonumber\\
&=-({\mathcal L}_{X}g)(\phi X',\phi X'')=0\nonumber,
\end{align}
\begin{align*}
(\nabla_{Y}g)(X',X'')&=Y(g(X',X''))-g([Y,X']_L,X'')-g(X',[Y,X'']_L)\\
&=Y(g(X',X''))-g([Y,X'],X'')-g(X',[Y,X''])\\
&=({\mathcal L}_{Y}g)(X',X'')=0
\end{align*}
and
\begin{align*}
(\nabla_{\xi}g)(X',X'')&=\xi(g(X',X''))-g([\xi,X']_L,X'')-g(X',[\xi,X'']_L)\\
&=\xi(g(X',X''))-g([\xi,X'],X'')-g(X',[\xi,X''])\\
&=({\mathcal L}_{\xi}g)(X',X'')=0
\end{align*}
for all $X,X',X''\in\Gamma(L)$ and $Y\in\Gamma(Q)$. Analogously one
can prove that $(\nabla_{V}g)(Y',Y'')=0$ for all
$V\in\Gamma(TM^{2n+1})$ and $Y',Y''\in\Gamma(Q)$.\\
Now we prove the last part of the theorem. We prove that, under the
assumption of the integrability of $L$ and $Q$, (i) is equivalent to
the total geodesicity of the foliations defined by $L$ and $Q$. Let
$X,X'$ be sections of $L$. Then for any $Y\in\Gamma(Q)$ the Koszul
formula for the Levi Civita connection yields
\begin{align}\label{equiv1}
2g(\hat\nabla_{X}X',Y)&=X(g(X',Y))+X'(g(X,Y))-Y(g(X,X'))+g([X,X'],Y)\nonumber\\
&\quad+g([Y,X],X')-g([X',Y],X)\nonumber\\
&=-Y(g(X,X'))+g([Y,X]_L,X')+g([Y,X']_L,X)\\
&=-Y(g(X,X'))+g(\nabla_{Y}X,X')+g(X,\nabla_{Y}X')\nonumber\\
&=-(\nabla_{Y}g)(X,X')\nonumber
\end{align}
and, in the same way,
\begin{equation}\label{equiv2}
2g(\hat\nabla_{X}X',\xi)=-(\nabla_{\xi}g)(X,X'),
\end{equation}
from which it follows that if the bi-Legendrian connection is metric
then the foliation defined by $L$ is totally geodesic. A similar
argument works also for $Q$. Conversely, if $L$ and $Q$ define two
totally geodesic foliations, by \eqref{equiv1}--\eqref{equiv2} one
has
$(\nabla_{Y}g)(X,X')=(\nabla_{X}g)(Y,Y')=(\nabla_{\xi}g)(X,X')=(\nabla_{\xi}g)(Y,Y')=0$
for any $X,X'\in\Gamma(L)$, $Y,Y'\in\Gamma(Q)$. Moreover, for all
$X,X',X''\in\Gamma(L)$, using the same computations in
\eqref{formulautile},
\begin{align*}
(\nabla_{X}g)(X',X'')&=-X(g(X',X''))-g(X',\phi[X,\phi X''])-g(X'',\phi[X,\phi X'])\\
&=\phi X'(g(\phi X'',X))+\phi X''(g(\phi X',X))-X(g(\phi X',\phi
X''))\\
&\quad +g([\phi X',\phi X''],X)+g([X,\phi X'],\phi X'')-g([\phi
X''X],\phi X')\\
&=2g(\hat\nabla_{\phi X'}\phi X'',X)=0
\end{align*}
because of the totally geodesicity of the foliation defined by $Q$.
Analogously one can prove that $(\nabla_{Y}g)(Y',Y'')=0$ for all
$Y,Y',Y''\in\Gamma(Q)$. Hence $\nabla g=0$.
\end{proof}

\begin{example}\label{kappamu}
\emph{A class of examples of bi-Legendrian structures verifying
one of the equivalent conditions stated in Proposition
\ref{metrica} is given by contact
$\left(\kappa,\mu\right)$-manifolds, i.e. contact metric manifolds
such that the Reeb vector field satisfies}
\begin{equation*}
\hat{R}(V,W)\xi=\kappa\left(\eta\left(W\right)V-\eta\left(V\right)W\right)+\mu\left(\eta\left(W\right)h
V-\eta\left(V\right)h W\right)
\end{equation*}
\emph{for some constants $\kappa, \mu \in\mathbb R$. This class of
contact metric manifolds has been introduced in  \cite{blair0} and
then extensively studied by several authors. It is well known
 that $\kappa\leq 1$ and when $\kappa<1$ the
contact metric manifold in question admits two mutually orthogonal
and integrable Legendrian distributions
${\mathcal{D}}\left(\lambda\right)$ and
${\mathcal{D}}\left(-\lambda\right)$ determined by the eigenspaces
of the operator $h$, where $\lambda=\sqrt{1-\kappa}$. Moreover,
these Legendrian foliations are totally geodesic, hence they verify
(i)--(iv) of Proposition \ref{metrica}. This bi-Legendrian structure
and the corresponding bi-Legendrian connection has been studied in
detail in \cite{mino3} where in particular it is proved that
${\mathcal{D}}\left(\lambda\right)$ and
${\mathcal{D}}\left(-\lambda\right)$ are never
 both flat.}
\end{example}

\section{The bi-Legendrian and the Tanaka-Webster connection}\label{confronto}

In this section we consider a contact metric manifold
$(M^{2n+1},\phi,\xi,\eta,g)$ endowed with a Legendrian
distribution $L$. We denote, as usual, by $Q$ the conjugate
Legendrian distribution of $L$ and by  $\nabla$ the bi-Legendrian
connection corresponding to $\left(L,Q\right)$. We assume that the
pair $\left(L,Q\right)$ is flat, that is both $L$ and $Q$ are flat
Legendrian distributions, and satisfies one of the equivalent four
properties of Proposition \ref{metrica}. Under these assumptions
we study the relationship between $\nabla$ and the Tanaka-Webster
connection $^\ast\nabla$ of
$(M^{2n+1},\phi,\xi,\eta,g)$.

\begin{theorem}\label{tanakabilegendrian}
Under the notation and the assumptions above, the bi-Legendrian
connection $\nabla$ coincides with the Tanaka-Webster connection
$^{\ast}\nabla$ if and only if $L$ and $Q$ are integrable and
$(M^{2n+1},\phi,\xi,\eta,g)$ is a Sasakian manifold.
\end{theorem}
\begin{proof}
Suppose that $\nabla={^{\ast}\nabla}$. Then the torsion tensor
field $T$ of the bi-Legendrian connection must satisfy (iv) of
\eqref{listatanaka}. In particular,
$\left[X,X'\right]_Q=-T\left(X,X'\right)=-2d\eta\left(X,X'\right){\xi}=0$
for all $X,X'\in\Gamma\left(L\right)$ and
$\left[Y,Y'\right]_L=-T\left(Y,Y'\right)=-2d\eta\left(Y,Y'\right){\xi}=0$
for all $Y,Y'\in\Gamma\left(Q\right)$, from which we deduce the
integrability of $L$ and $Q$. Moreover, from \eqref{tanakatorsion}
it follows that $\left[\xi,X\right]_Q=T\left(X,\xi\right)={^\ast
T}\left(X,\xi\right)=\eta\left(\xi\right)\phi h
X-\eta\left(X\right)\phi h \xi+2g\left(X,\phi\xi\right)\xi=-h\phi
X$. So, for all $X\in\Gamma\left(L\right)$,
\begin{equation}\label{acca1}
\left[\xi,X\right]_{Q}=-h\phi X
\end{equation}
and, in the same way,
\begin{equation}\label{acca2}
\left[\xi,Y\right]_{L}=-h\phi Y
\end{equation}
for all $Y\in\Gamma\left(Q\right)$. By \eqref{acca1} and
\eqref{acca2} we see that the flatness of $L$ and $Q$ is equivalent
to the vanishing of $h$. With this remark we can prove that
$(M^{2n+1},\phi,\xi,\eta,g)$ is  Sasakian. Indeed, since $\nabla
g=0$, by Proposition \ref{metrica} we have $\nabla\phi=0$. Moreover,
$\nabla$ satisfies (ii) of \eqref{listatanaka}, so for all
$V,W\in\Gamma\left(TM\right)$
\begin{equation*}
(\hat{\nabla}_V\phi)W=g\left(V+hV,W\right)\xi-\eta\left(W\right)\left(V+hV\right)=g\left(V,W\right)\xi-\eta\left(W\right)V,
\end{equation*}
since $h=0$. Now we prove the converse, showing that $\nabla$
verifies \eqref{listatanaka}. We already know that $\nabla$
satisfies $\nabla\xi=0$, $\nabla\eta=0$ and, by hypothesis, $\nabla
g=0$. Moreover $\nabla$ satisfies also
$T\left(X,Y\right)=2d\eta\left(X,Y\right){\xi}$ for all
$X\in\Gamma\left(L\right)$ and $Y\in\Gamma\left(Q\right)$, so in
order to check (iv) it is sufficient to prove that
$T\left(X,X'\right)=T\left(Y,Y'\right)=0$ for all
$X,X'\in\Gamma\left(L\right)$, $Y,Y'\in\Gamma\left(Q\right)$. But
this is true because, by the assumption of the integrability of $L$
and $Q$, we have $T\left(X,X'\right)=-\left[X,X'\right]_{Q}=0$ and
$T\left(Y,Y'\right)=-\left[Y,Y'\right]_{L}=0$. Moreover, since
$(M^{2n+1},\phi,\xi,\eta,g)$ is a Sasakian manifold and, in
particular, a K-contact manifold, we have, for all
$V,W\in\Gamma\left(TM\right)$,
\begin{gather*}
(\hat{\nabla}_V\phi)W-g\left(V+hV,W\right)\xi+\eta\left(W\right)\left(V+hV\right)\\
=(\hat{\nabla}_V\phi)W-g\left(V,W\right)\xi+\eta\left(W\right)\left(V\right)=0=\left(\nabla_{V}\phi\right)W
\end{gather*}
because of Proposition \ref{metrica}. So $\nabla$ satisfies also
(ii). Finally, since $h=0$ and $L,Q$ are flat, we have, for all
$X\in\Gamma\left(L\right)$, $T\left(\xi,\phi X\right)=\left[\phi
X,\xi\right]_L=0=-\phi(\left[X,\xi\right]_Q)=-\phi
T\left(\xi,X\right)$, and, similarly, for all
$Y\in\Gamma\left(Q\right)$, $T\left(\xi,\phi Y\right)=0=-\phi
T\left(\xi,Y\right)$, hence also (iii) is satisfied. Thus by the
uniqueness of the Tanaka-Webster connection, we conclude that
$\nabla={^{\ast}\nabla}$.
\end{proof}

\begin{remark}
\emph{In the proof of Theorem \ref{tanakabilegendrian}, we have
found the following expression for the tensor field $h$:}
\begin{equation*}
h X = \left[\xi,\phi
X\right]_L=-\left(\phi\left[\xi,X\right]\right)_L, \ h Y =
\left[\xi,\phi Y\right]_Q=-\left(\phi\left[\xi,Y\right]\right)_Q,
\end{equation*}
\emph{for all $X\in\Gamma\left(L\right)$ and
$Y\in\Gamma\left(Q\right)$. In particular, as we already know by
Proposition \ref{metrica}, $h$ preserves the distributions $L$ and
$Q$.}
\end{remark}

As immediate consequences of Theorem \ref{tanakabilegendrian} and
Proposition \ref{metrica} we have:

\begin{corollary}
Under the assumptions of Theorem \ref{tanakabilegendrian}, the
Tanaka-Webster connection of $(M^{2n+1},\phi,\xi,\eta,g)$ satisfies
$^\ast\nabla_{X}X'=-\left(\phi\left[X,\phi X'\right]\right)_L$ for
all $X,X'$ $\in\Gamma\left(L\right)$ and
$^\ast\nabla_{Y}Y'=-\left(\phi\left[Y,\phi Y'\right]\right)_Q$ for
all $Y,Y'\in\Gamma\left(Q\right)$.
\end{corollary}

\begin{corollary}\label{levicivita}
Let $(M^{2n+1},\phi,\xi,\eta,g)$ be a Sasakian manifold foliated by
a flat Legendrian foliation $\cal{F}$ such that the conjugate
Legendrian distribution is integrable. Let $\nabla$ be the
corresponding bi-Legendrian connection and suppose that $\nabla
g=0$. Let $S$ be the tensor field of type $(1,2)$ defined by
$S\left(V,W\right)=\nabla_{V}W-\hat{\nabla}_{V}W$. Then we have
$S\left(V,\xi\right)=S\left(\xi,V\right)=\phi V$
 for all  $V\in\Gamma\left(TM\right)$ and
$S\left(Z,Z'\right)=d\eta\left(Z,Z'\right){\xi}$  for all
$Z,Z'\in\Gamma\left(\cal{D}\right)$. In particular, for all
$X,X'\in\Gamma\left(L\right)$ and for all
$Y,Y'\in\Gamma\left(Q\right)$ we have
\begin{equation}\label{bilegendrianlevicivita}
\nabla_{X}X'=\hat{\nabla}_{X}X', \ \
\nabla_{Y}Y'=\hat{\nabla}_{Y}Y'.
\end{equation}
\end{corollary}
\begin{proof}
Indeed, by Theorem \ref{tanakabilegendrian}, $\nabla$ coincides
with the Tanaka-Webster connection of
$(M^{2n+1},\phi,\xi,\eta,g)$. Then, by
\eqref{tanakaconnection} we deduce the following relations:
\begin{gather*}
\nabla_{X}X'-\hat{\nabla}_{X}X'=0 \textrm{,  } \
\nabla_{X}Y'-\hat{\nabla}_{X}Y=d\eta\left(X,Y\right){\xi}
\textrm{,  } \
\nabla_{X}\xi-\hat{\nabla}_{X}\xi=\phi X,\\
\nabla_{Y}X-\hat{\nabla}_{Y}X=d\eta\left(Y,X\right){\xi} \textrm{,
} \ \nabla_{Y}Y'-\hat{\nabla}_{Y}Y=0 \textrm{,  } \
\nabla_{Y}\xi-\hat{\nabla}_{Y}\xi=\phi Y,\\
\nabla_{\xi}X-\hat{\nabla}_{\xi}X=\phi X \textrm{, } \
\nabla_{\xi}Y-\hat{\nabla}_{\xi}Y=\phi Y \textrm{, } \
\nabla_{\xi}\xi=\hat{\nabla}_{\xi}\xi=0
\end{gather*}
for all $X,X'\in\Gamma\left(L\right)$,
$Y,Y'\in\Gamma\left(Q\right)$, from which the assertion follows.
\end{proof}

\begin{remark}\label{osservazione}
\emph{Let $(M^{2n+1},\phi,\xi,\eta,g)$ be a Sasakian manifold and
let $\Im_M$ be the set of all flat Legendrian foliations on
$M^{2n+1}$ such that the conjugate Legendrian distribution is
integrable and $\nabla g=0$. Take two elements ${\cal{F}}_1$ and
${\cal{F}}_2$ of $\Im_M$. ${\cal{F}}_1$ and ${\cal{F}}_2$ are flat
Legendrian foliations on $M^{2n+1}$ such that $\nabla^1 g=0$ and
$\nabla^2 g=0$, where $\nabla^1$ and $\nabla^2$ denote the
bi-Legendrian connections associated to ${\cal{F}}_1$ and
${\cal{F}}_2$, respectively. Then, by Theorem
\ref{tanakabilegendrian}, $\nabla^1=\nabla^2$ because they both
coincide with the Tanaka-Webster connection $^{\ast}\nabla$. In
particular we have that $\nabla^1{\cal{F}}_2\subset{\cal{F}}_2$ and
$\nabla^2{\cal{F}}_1\subset{\cal{F}}_1$. Moreover we deduce that the
Tanaka-Webster connection preserves all the Legendrian foliations
belonging to $\Im_M$.}
\end{remark}

A variation of Theorem \ref{tanakabilegendrian} is the following
Theorem \ref{tanakabilegendrian1}. But, before proving it, we need a
preliminary lemma.

\begin{lemma}\label{flat}
Let $(M^{2n+1},\phi,\xi,\eta,g)$ be a $K$-contact manifold endowed
with a flat Legendrian distribution $L$. Then its conjugate
Legendrian distribution $Q=\phi L$ is also flat.
\end{lemma}
\begin{proof}
Indeed, as $\xi$ is Killing, we have $h=0$, so that, for all
$X\in\Gamma\left(L\right)$, $0=2hX=\left[\xi,\phi
X\right]-\phi\left[\xi,X\right]$, from which $\left[\xi,\phi
X\right]=\phi\left[\xi,X\right]\in\Gamma\left(Q\right)$, because
$L$ is flat.
\end{proof}

\begin{theorem}\label{tanakabilegendrian1}
Let $(M^{2n+1},\phi,\xi,\eta,g)$ be a Sasakian manifold endowed with
a flat Legendrian distribution $L$. Let $Q=\phi L$ be its conjugate
Legendrian distribution. If the Tanaka-Webster connection
$^\ast\nabla$ preserves the distribution $L$, then $L$ and $Q$ are
integrable and $^\ast\nabla$ coincides with the bi-Legendrian
connection $\nabla$ associated to the almost bi-Legendrian structure
$\left(L,Q\right)$.
\end{theorem}
\begin{proof}
First of all, we prove that $^\ast\nabla L\subset L$ implies the
integrability of $L$. Let $X,X'\in\Gamma\left(L\right)$. Then
$\left[X,X'\right]={^\ast\nabla}_{X}X'-{^\ast\nabla}_{X'}X-2d\eta\left(X,X'\right)\xi={^\ast\nabla}_{X}X'-{^\ast\nabla}_{X'}X\in\Gamma\left(L\right)$.
Now we show that $^\ast\nabla Q\subset Q$. Let
$Y\in\Gamma\left(Q\right)$. Then, by $^\ast\nabla g=0$ and
${^\ast\nabla} L\subset L$, we get for all $V\in\Gamma(TM^{2n+1})$
and $X\in\Gamma\left(L\right)$
\begin{align*}
0=\left(^\ast\nabla_{V}g\right)\left(X,Y\right)=V\left(g\left(X,Y\right)\right)-g\left(^\ast\nabla_{V}X,Y\right)-g\left(X,^\ast\nabla_{V}Y\right)=-g\left(X,^\ast\nabla_{V}Y\right),
\end{align*}
so that $^\ast\nabla_{V}Y\in\Gamma(Q\oplus\mathbb{R}\xi)$. Moreover,
since $^\ast\nabla\xi=0$,
$0=\left(^\ast\nabla_{V}g\right)\left(\xi,Y\right)=V\left(g\left(\xi,Y\right)\right)-g\left(^\ast\nabla_{V}\xi,Y\right)-g\left(\xi,^\ast\nabla_{V}Y\right)=-g\left(\xi,^\ast\nabla_{V}Y\right)$,
from which $^\ast\nabla_{V}Y\in\Gamma\left(Q\right)$. Then, arguing
in the same way as for $L$, one can prove that $Q$ is integrable.
Note also that since $M^{2n+1}$ is Sasakian and in particular
K-contact, by Lemma \ref{flat}, also $Q$ is flat. Finally, we prove
that $^\ast\nabla$ coincides with the bi-Legendrian connection
corresponding to $\left(L,Q\right)$, that is $^\ast\nabla$ verifies
(ii) and (iii) in \eqref{lista}. The relations  $^\ast
T\left(X,\xi\right)=\left[\xi,X\right]_Q$ for
$X\in\Gamma\left(L\right)$ and $^\ast
T\left(Y,\xi\right)=\left[\xi,Y\right]_L$ for
$Y\in\Gamma\left(Q\right)$ hold because $L$ and $Q$ are flat and, on
the other hand, $^\ast T\left(X,\xi\right)=\phi h X=0$, $^\ast
T\left(Y,\xi\right)=\phi h Y=0$. In order to prove $^\ast\nabla
d\eta=0$, we show firstly that $^\ast\nabla\phi=0$. Indeed, since
$M^{2n+1}$ is Sasakian,
\begin{equation*}
\left(^\ast\nabla_{V}\phi\right)W=(\hat{\nabla}_{V}\phi)W-g\left(V,W\right)\xi+\eta\left(W\right)V=0
\end{equation*}
for all $V,W\in\Gamma(TM^{2n+1})$, so that $^\ast\nabla\phi=0$. Now
we can prove that $\left(^\ast\nabla_{V}
d\eta\right)\left(W,W'\right)$ $=0$ for all
$V,W,W'\in\Gamma(TM^{2n+1})$. This equality holds immediately for
$W,W'\in\Gamma\left(L\right)$ and for $W,W'\in\Gamma\left(Q\right)$
because $L$ and $Q$ are preserved by $^\ast\nabla$. Also the case
$W'=\xi$ is obvious since $^\ast\nabla\xi=0$. So it remains to show
that $\left(^\ast\nabla_{V} d\eta\right)\left(X,Y\right)=0$ for
$X\in\Gamma\left(L\right)$ and $Y\in\Gamma\left(Q\right)$. In fact,
using $^\ast\nabla\phi=0$,
\begin{align*}
\left(^\ast\nabla_{V} d\eta\right)\left(X,Y\right)&=V\left(g\left(X,\phi Y\right)\right)-g\left(^\ast\nabla_{V}X,\phi Y\right)-g\left(X,\phi{^\ast\nabla}_{V}Y\right)\\
&=V\left(g\left(X,\phi Y\right)\right)-g\left(^\ast\nabla_{V}X,\phi Y\right)-g\left(X,^\ast\nabla_{V}\phi Y\right)\\
&=\left(^\ast\nabla_{V}g\right)\left(X,\phi Y\right)=0,
\end{align*}
since $^\ast\nabla g=0$. Thus $^\ast\nabla$ satisfies all the
properties which characterize the bi-Legendrian connection
associated to $\left(L,Q\right)$.
\end{proof}

\section{An interpretation of the Tanaka-Webster connection}\label{interpretazione}

In $\S$ \ref{confronto} we have found that under certain assumptions
the Tanaka-Webster connection of a Sasakian manifold foliated by a
Legendrian foliation $\cal F$ coincides with the bi-Legendrian
connection associated to $\cal F$ (Theorem
\ref{tanakabilegendrian1}). This result has an analogue in even
dimension: F. Etayo and R. Santamaria proved in \cite{etayo1} that
under suitable assumptions the Levi Civita connection of a
K\"{a}hlerian manifold foliated by a Lagrangian foliation $\cal F'$
coincides with the bi-Lagrangian connection associated to $\cal F'$.
Therefore it seems that the Tanaka-Webster connection plays the same
role of the Levi Civita connection for symplectic or K\"{a}hlerian
manifolds. This is not surprising, since it is a well-known fact
that the Tanaka-Webster connection of a Sasakian manifold which is a
circle bundle over a K\"{a}hlerian manifold can be viewed as the
lift of the Levi Civita connection of the K\"{a}hlerian base
manifold. Now we prove this property for any, in general
non-regular, Sasakian manifold.

Let $(M^{2n+1},\phi,\xi,\eta,g)$ be a Sasakian manifold. It is well
known that the Reeb vector field $\xi$ defines a transversely
K\"{a}hlerian foliation, that is this foliation, which we denote by
${\cal F}_{\xi}$, can be defined by local submersions
$f_i:U_i\longrightarrow M'^{2n}$ from an open set $U_i$ of
$M^{2n+1}$, with $\left\{U_i\right\}_{i\in I}$ an open covering of
$M^{2n+1}$, onto a K\"{a}hlerian manifold $(M'^{2n},J,\omega,G)$,
where $J$, $\omega$ and $G$ are the projection of $\phi$, $d\eta$
and $g$, respectively. Moreover, any two of these submersions $f_i$
and $f_j$, with $U_i\cap U_j\neq\emptyset$, are connected by local
diffeomorphisms $\gamma_{ij}:f_j(U_i\cap U_j)\longrightarrow
f_i(U_i\cap U_j)$ satisfying, on $U_i\cap U_j$, the relation
$\gamma_{ij}\circ f_j=f_i$, and preserving the K\"{a}hlerian
structure.  Let $\hat{\nabla}'$ be the Levi Civita connection of
$(M'^{2n},G)$ and define a connection $\nabla^i$, locally on each
$U_i$, as the lift of $\hat{\nabla}'$ under the submersion $f_i$.
More precisely, for any basic vector fields $Z_1,Z_2$, we define
$\nabla^i_{Z_1}Z_2$ as the unique basic vector field on $U_i$ such
that $f_{i \ast}(\nabla^i_{Z_1}Z_2)=\hat{\nabla}'_{f_{i
\ast}(Z_1)}f_{i \ast}(Z_2)$. Moreover, we put, by definition,
$\nabla^i\xi=0$ and, for any vector field $V$ on $U_i$,
$\nabla^i_{\xi}V=\left[\xi,V\right]$. Note that these last
definitions implies that, for any basic vector field $Z$, $f_{i
\ast}(\nabla^i_{Z}\xi)=0=\hat{\nabla}'_{f_{i
\ast}\left(Z\right)}f_{i \ast}\left(\xi\right)$ and $f_{i
\ast}(\nabla^i_{\xi}Z)=f_{i
\ast}\left(\left[\xi,Z\right]\right)=0=\hat{\nabla}'_{f_{i
\ast}\left(\xi\right)}f_{i \ast}\left(Z\right)$. Note also that
$\nabla^i$ preserves the "horizontal" distribution $\cal D$. We have
the following result:

\begin{proposition}\label{sollevamento}
The above connection $\nabla^i$ coincides with the Tanaka-Webster
connection of $(M^{2n+1},\phi,\xi,\eta,g)$ restricted
to $U_i$.
\end{proposition}
\begin{proof}
It is sufficient to show that $\nabla^i$ verifies all the
properties which characterize the Tanaka-Webster connection of
$(M^{2n+1},\phi,\xi,\eta,g)$. First of all, by
definition, $\nabla^i \xi=0$. Next, from $\nabla^i\xi=0$ and
$\nabla^i{\cal{D}}\subset{\cal{D}}$, we deduce $\nabla^i\eta=0$.
Furthermore, since $\hat{\nabla}'G=0$ and $f_i$ is a Riemannian
submersion, we get
$(\nabla^i_{Z}g)(Z_1,Z_2)=0$ for all
$Z,Z_1,Z_2$ basic vector fields on $U_i$, and, since $\nabla^i\cal
D \subset \cal D$, also
$(\nabla^i_{Z}g)(Z_1,\xi)=0$. So it remains
to prove that $(\nabla^i_{\xi}g)(Z_1,Z_2)=0$ for
$Z_1,Z_2$ basic vector fields. Indeed
\begin{equation*}
(\nabla^i_{\xi}g)(Z_1,Z_2)=\xi\left(g(Z_1,Z_2)\right)-g\left(\left[{\xi},Z_1\right],Z_2\right)-g\left(Z_1,\left[{\xi},Z_2\right]\right)=\left({\cal{L}}_{\xi}g\right)\left(Z_1,Z_2\right)=0
\end{equation*}
because $\xi$ is Killing. In the same way, since $\hat{\nabla}' J=0$
and $f_{i \ast}\circ \phi=J\circ f_{i \ast}$, we get
$(\nabla^{i}_{Z_1}\phi)Z_2=0$ for all $Z_1,Z_2$ basic vector fields
on $U_i$. Next, for any basic vector field $Z$ on $U_i$ we have
$(\nabla^{i}_{\xi}\phi)Z=\left[\xi,\phi
Z\right]-\phi\left[\xi,Z\right]=2h Z=0$ because $h=0$, $M^{2n+1}$
being Sasakian. Thus for concluding the proof it remains to check
the properties involving the torsion. Let $Z$ be a basic vector
field defined on $U_i$. Then $T^i\left(\xi,\phi
Z\right)=\nabla^i_{\xi}\phi Z-\nabla^i_{\phi Z}\xi-\left[\xi,\phi
Z\right]=\left[\xi,\phi Z\right]-\left[\xi,\phi Z\right]=0$ and
$T^i\left(\xi, Z\right)=\nabla^i_{\xi} Z-\nabla^i_{ Z}\xi-\left[\xi,
Z\right]=\left[\xi, Z\right]-\left[\xi, Z\right]=0$, so that
$T^i\left(\xi,\phi Z\right)=0=-\phi T^i\left(\xi, Z\right)$.
Finally, for any $Z_1,Z_2$ basic vector fields, we have $f_{i
\ast}(T^i(Z_1,Z_2))=T'\left(f_{i \ast}\left(Z_1\right),f_{i
\ast}\left(Z_2\right)\right)=0$ and so $T^i\left(Z_1,Z_2\right)$ is
vertical. Hence
$T^i\left(Z_1,Z_2\right)=\eta\left(T^i\left(Z_1,Z_2\right)\right)\xi=-\eta\left(\left[Z_1,Z_2\right]\right)\xi=2d\eta\left(Z_1,Z_2\right)\xi$.
\end{proof}
\\

Now we prove that this family of connections give rise to a
well-defined global connection on $M^{2n+1}$.

\begin{proposition}\label{globale}
Let $(M^{2n+1},\phi,\xi,\eta,g)$ be a Sasakian manifold and let
$\left\{U_i,f_i,\gamma_{ij}\right\}$ be a cocycle defining the
foliation ${\cal F}_\xi$. Then the family of connections
$(\nabla^i)_{i\in I}$ defined above gives rise to a global
connection on $M^{2n+1}$ which coincides with the Tanaka-Webster
connection of $(M^{2n+1},\phi,\xi,\eta,g)$.
\end{proposition}
\begin{proof}
We have to prove that for any $i,j\in I$ such that $U_i\cap
U_j\neq\emptyset$, $\nabla^i=\nabla^j$. Firstly note that
$\gamma_{ij}:f_j(U_i\cap U_j)\longrightarrow f_i(U_i\cap U_j)$ is a
(local) affine transformation with respect to the Levi Civita
connection, because it is a (local) isometry.  Now let $Z_1'$ and
$Z_2'$ be vector fields on $M'^{2n}$ and let $Z_1^i$, $Z_1^j$ and
$Z_2^i$, $Z_2^j$ be the basic vector fields $f_i$-related and
$f_j$-related, respectively, to $Z_1'$ and $Z_2'$. Note that
$Z_1^{i}$ is also the basic vector field $f_{j}$-related to
$\gamma_{ij_\ast}(Z_1^{\prime})$ because it is horizontal also for
$f_j$, as
$\ker\left(f_{i\ast}\right)=\mathbb{R}\xi=\ker\left(f_{j\ast}\right)$,
and, for all $p\in M^{2n+1}$,
$f_{j_\ast}(Z_{1_p}^{i})=(\gamma_{ij_\ast})_{f_j\left(p\right)}(Z'_{1_{f_j\left(p\right)}})$,
since $f_i\left(p\right)=\gamma_{ij}\left(f_j\left(p\right)\right)$
and $\gamma_{ij}=\gamma_{ji}^{-1}$. Then we get
$f_{i_\ast}(\nabla^{i}_{Z^i_{1}}Z^i_{2})=\gamma_{ij_\ast}(f_{j_\ast}(\nabla^j_{Z^i_1}Z^i_2))=f_{i_\ast}(\nabla^{j}_{Z^i_{1}}Z^i_{2})$,
which implies that $\nabla^i_{Z^i_1}Z^i_2-\nabla^j_{Z^i_1}Z^i_2$ is
vertical. Since it is also horizontal, we get
$\nabla^{i}_{Z^i_{1}}Z^i_{2}=\nabla^{j}_{Z^i_{1}}Z^i_{2}$. Moreover,
clearly, $\nabla^i\xi=0=\nabla^j\xi$ and, on $U_i \cap U_j$,
$\nabla^i_{\xi}V=\left[\xi,V\right]=\nabla^j_{\xi}V$. Finally, the
last part of the statement follows from Proposition
\ref{sollevamento}.
\end{proof}
\\

More in general, for any contact metric manifolds
$(M^{2n+1},\phi,\xi,\eta,g)$ we can define a connection
on $M^{2n+1}$ setting, for all $Z\in\Gamma\left(\cal D\right)$,
\begin{equation}\label{connessionecanonica}
\tilde{\nabla}_{Z}Z'=(\hat{\nabla}_ZZ')_{\cal D}, \
\tilde{\nabla}_{\xi}Z=\left[\xi,Z\right], \ \tilde{\nabla}\xi=0.
\end{equation}
That $\tilde{\nabla}$ is a connection on $M^{2n+1}$ preserving the contact
distribution $\cal D$ is easy to check. Moreover, we can give an
interesting characterization of this connection:

\begin{theorem}
Let $(M^{2n+1},\phi,\xi,\eta,g)$ be a contact metric manifold and
$\tilde{\nabla}$ the connection on $M^{2n+1}$ defined by
\eqref{connessionecanonica}. Then $\tilde{\nabla}$ is the unique
connection on $M^{2n+1}$ satisfying the following properties:
\begin{description}
    \item[(i)] $\tilde{\nabla}\xi=0$,
    \item[(ii)]
    $\tilde{T}\left(V,W\right)=2\eta\left(V,W\right)\xi$ for all
    $V,W\in\Gamma\left(TM\right)$,
    \item[(iii)] $(\tilde{\nabla}_{Z}g)\left(Z',Z''\right)=0$
    for all $Z,Z',Z''\in\Gamma\left(\cal D\right)$.
\end{description}
Furthermore, $M^{2n+1}$ is K-contact if and only if
$\tilde{\nabla}g=0$, and $M^{2n+1}$ is Sasakian if and only if
$\tilde{\nabla}\phi=0$ and in this case $\tilde{\nabla}$ coincides
with the Tanaka-Webster connection of $(M^{2n+1},\phi,\xi,\eta,g)$.
\end{theorem}
\begin{proof}
Firstly we prove that $\tilde{\nabla}$ satisfies (i), (ii) and
(iii). By definition we have $\tilde{\nabla}\xi=0$. Next, for all
$Z,Z'\in\Gamma\left(\cal D\right)$,
$\tilde{T}\left(Z,Z'\right)=(\hat{\nabla}_{Z}Z')_{\cal
D}-(\hat{\nabla}_{Z'}Z)_{\cal
D}-\left[Z,Z'\right]=(\hat{T}\left(Z,Z'\right))_{\cal
D}-\left[Z,Z'\right]_{\mathbb{R}\xi}=-\eta\left(\left[Z,Z'\right]\right)\xi=2d\eta\left(Z,Z'\right)\xi$,
and
$\tilde{T}\left(Z,\xi\right)=-\tilde{\nabla}_{\xi}Z-\left[Z,\xi\right]=\left[Z,\xi\right]-\left[Z,\xi\right]=0=2d\eta\left(Z,\xi\right)\xi$.
Then, since on the contact distribution $\tilde{\nabla}$ coincides with the projection on $\cal D$ of the Levi
Civita connection, we get (iii). Now let $\nabla$ be any connection
on $M^{2n+1}$ satisfying (i), (ii), (iii). Then by (i) and (ii) we
have $\nabla_\xi
Z=\nabla_Z\xi+\left[\xi,Z\right]+T\left(\xi,Z\right)=\left[\xi,Z\right]+2d\eta\left(\xi,Z\right)\xi=\left[\xi,Z\right]$
for all $Z\in\Gamma\left(\cal D\right)$. So it remains to prove
that $\nabla_Z Z'=(\hat{\nabla}_Z Z')_{\cal D}$ for all
$Z,Z'\in\Gamma\left(\cal D\right)$. For this purpose, let
$\bar{\nabla}$ be the connection given by
\begin{equation*}
\bar{\nabla}_{V}W=\nabla_{V_{\cal D}}W_{\cal
D}+(\hat{\nabla}_{V_{\cal D}}W_{\cal
D})_{\mathbb{R}\xi}+\hat{\nabla}_{V_{\mathbb{R}\xi}}W+\hat{\nabla}_{V}{W_{\mathbb{R}\xi}}.
\end{equation*}
Then, if we prove that $\bar{\nabla}$ coincides with the Levi
Civita connection of $M^{2n+1}$, we would have that $\nabla_Z
Z'=(\hat{\nabla}_Z Z')_{\cal D}$ for all $Z,Z'\in\Gamma\left(\cal
D\right)$. It is enough to verify that $\bar{\nabla}$ is
metric and torsion free on the subbundle $\cal D$. That
$\bar{\nabla}$ is metric on $\cal D$ is ensured by (iii); then,
$\bar{T}\left(Z,Z'\right)=T\left(Z,Z'\right)+\eta(\hat{\nabla}_{Z}Z')\xi-\eta(\hat{\nabla}_{Z'}Z)\xi=2d\eta\left(Z,Z'\right)\xi+\eta\left(\left[Z,Z'\right]\right)\xi=0$.
For proving the second part of the theorem, note that
\begin{equation*}
(\tilde{\nabla}_{\xi}g)(Z,Z')=\xi(g(Z,Z'))-g([\xi,Z],Z')-g(Z,[\xi,Z'])=({\cal
L}_{\xi}g)(Z,Z'),
\end{equation*}
from which we deduce that $M^{2n+1}$ is a K-contact manifold if
and only if $\tilde{\nabla}$ is a metric connection with respect to
the associated metric $g$. Finally, if $\tilde{\nabla}\phi=0$ we have, first of all,
\begin{equation}\label{kcontatto}
0=(\tilde{\nabla}_{\xi}\phi)Z=\left[\xi,\phi
Z\right]-\phi\left[\xi,Z\right]=\left({\cal
L}_{\xi}\phi\right)Z=2hZ,
\end{equation}
from which $M^{2n+1}$ is K-contact and by \eqref{condizione} $\phi
Z=-\hat{\nabla}_Z\xi$. Then, for all $Z,Z'\in\Gamma\left(\cal
D\right)$, $(\hat{\nabla}_Z\phi)Z'=((\hat{\nabla}_Z\phi)Z')_{\cal
D}+((\hat{\nabla}_Z\phi)Z')_{\mathbb{R}\xi}=(\tilde{{\nabla}}_Z\phi)Z'+\eta((\hat{\nabla}_Z\phi)Z')\xi=g(\hat{\nabla}_{Z}\phi
Z',\xi)\xi=-g(\phi Z',\hat{\nabla}_Z\xi)\xi=g\left(\phi Z,\phi
Z'\right)\xi=g\left(Z,Z'\right)\xi$, and \eqref{condizionesasaki} is
satisfied. Moreover, $(\hat{\nabla}_{\xi}\phi)Z=\hat{\nabla}_{\phi
Z}\xi+\left[\xi,\phi
Z\right]-\phi\hat{\nabla}_Z\xi-\phi\left[\xi,Z\right]=-\phi^2Z+\phi^2Z+\left({\cal
L}_{\xi}\phi\right)Z=0$ and
$(\hat{\nabla}_{Z}\phi)\xi=-\phi\hat{\nabla}_Z\xi=\phi^2 Z=-Z$, so
that \eqref{condizionesasaki} holds in any case. Conversely, if
$M^{2n+1}$ is Sasakian, then it is K-contact hence
$(\tilde{\nabla}_{\xi}\phi)Z=({\cal L}_{\xi}\phi)Z=0$; moreover, for
any $Z,Z'\in\Gamma\left(\cal D\right)$,
$(\tilde{\nabla}_{Z}\phi)Z'=(g\left(Z,Z'\right)\xi-\eta\left(Z'\right)Z)_{\cal
D}=0$. Finally, Proposition \ref{globale} implies that
$\tilde{\nabla}$ is the Tanaka-Webster connection of the Sasakian
manifold $(M^{2n+1},\phi,\xi,\eta,g)$.
\end{proof}

\bigskip

In the context of symplectic geometry, in the Appendix we shall
prove the following result.

\begin{theorem}\label{kahler1}
Let $(M^{2n},\omega)$ be a symplectic manifold endowed
with a bi-Lagrangian structure $\left(\cal{F},\cal{G}\right)$ such
that $T\cal{G}$ is an affine transversal distribution for $\cal
F$. Then there exists a K\"{a}hlerian structure on
$(M^{2n},\omega)$ whose Levi Civita connection
coincides with the bi-Lagrangian connection of
$(M^{2n},\omega,\cal{F},\cal{G})$.
\end{theorem}

Now we prove the analogue of Theorem \ref{kahler1} in odd
dimension. As it is expected, the role played in Theorem
\ref{kahler1} by the Levi Civita connection is played now by the
Tanaka-Webster connection:

\begin{theorem}\label{sasaki1}
Let $(M^{2n+1},\eta)$ be a contact manifold endowed with a flat
bi-Legendrian structure $\left(\cal{F},\cal{G}\right)$ such that
$T\cal{G}$ is an affine transversal distribution for $\cal F$. Then
there exists a Sasakian structure on $(M^{2n+1},\eta)$ whose
Tanaka-Webster connection coincides with the bi-Legendrian
connection of $(M^{2n+1},\eta,\cal{F},\cal{G})$.
\end{theorem}
\begin{proof}
The assumption of $T\cal G$ being an affine transversal distribution
for $\cal F$ means that the curvature tensor field of the
corresponding bi-Legendrian connection satisfies
$R\left(X,Y\right)=0$ for $X\in\Gamma\left(T\cal F\right)$,
$Y\in\Gamma\left(T\cal G\right)$ (cf. \cite{cappellettimontano1}).
So this assumption and the flatness of the bi-Lengendrian structure
imply that the curvature $R$ of the bi-Legendrian connection
$\nabla$ associated to $\left(\cal F, \cal G\right)$ vanishes
identically (cf. Proposition \ref{proprieta}). Now  let $p$ be a
point of $M^{2n+1}$. Since $d\eta_{p}$ is a symplectic form on the
subspace ${\cal{D}}_{p}\subset T_{p}M^{2n+1}$, it follows that there
exists a basis
$\left\{e_{1},\ldots,e_{n},e_{n+1},\ldots,e_{2n},\xi_{p}\right\}$ of
$T_{p}M^{2n+1}$ such that $\left\{e_{1},\ldots,e_{n}\right\}$ is a
basis of $T_{p}\cal F$, $\left\{e_{n+1},\ldots,e_{2n}\right\}$ is a
basis of $T_{p}\cal G$ and
\begin{equation}\label{relazioni}
d\eta_p\left(e_{i},e_{j}\right)=d\eta_p\left(e_{n+i},e_{n+j}\right)=0,  \
\ d\eta_p\left(e_{i},e_{n+j}\right)=-\frac{1}{2}\delta_{ij}
\end{equation}
for all $i,j\in\left\{1,\ldots,n\right\}$. For each
$k\in\left\{1,\ldots,2n\right\}$ we define vector fields $E_k$ on
$M^{2n+1}$ by the $\nabla$-parallel transport along curves of the
vector $e_k$. More precisely, for any $q\in M^{2n+1}$ we consider a
curve $\gamma:\left[0,1\right]\longrightarrow M$ such that
$\gamma\left(0\right)=p$, $\gamma\left(1\right)=q$ and we define
$E_k\left(q\right):=\tau_\gamma\left(e_k\right)$,
$\tau_\gamma:T_pM^{2n+1}\longrightarrow T_qM^{2n+1}$ being the
parallel transport along $\gamma$. Note that $E_{k}\left(q\right)$
does not depend on the curve joining $p$ and $q$, since $R\equiv 0$.
Setting $X_i:=E_{n+i}$ and $Y_i:=E_i$, we obtain $2n$ vector fields
on $M^{2n+1}$ such that, for each $i\in\left\{1,\ldots,n\right\}$,
$Y_{i}\in\Gamma\left(T\cal F\right)$ and $X_{i}\in\Gamma\left(T\cal
G\right)$, since the parallel transport preserves the distributions
$T\cal{F}$ and $T\cal{G}$. Moreover, \eqref{relazioni} holds at any
point of $M^{2n+1}$, that is for any $q\in M^{2n+1}$ and
$i,j\in\left\{1,\ldots,n\right\}$
\begin{equation}\label{relazioni2}
d\eta_q\left(Y_{i}\left(q\right),Y_{j}\left(q\right)\right)=d\eta_q\left(X_{i}\left(q\right),X_{j}\left(q\right)\right)=0, \ \
d\eta_q\left(Y_{i}\left(q\right),X_{j}\left(q\right)\right)=-\frac{1}{2}\delta_{ij}
\end{equation}
Indeed, since $d\eta$ is parallel with respect to $\nabla$, for
all $h,k\in\left\{1,\ldots,2n\right\}$,
\begin{equation*}
\frac{d}{dt}d\eta_{\gamma\left(t\right)}\left(E_{h}\left(\gamma\left(t\right)\right),E_{k}\left(\gamma\left(t\right)\right)\right)=d\eta_{\gamma\left(t\right)}\left(\nabla_{\gamma'}E_{h},E_{k}\right)+d\eta_{\gamma\left(t\right)}\left(E_{h},\nabla_{\gamma'}E_{k}\right)=0
\end{equation*}
so that
$d\eta_{p}\left(e_{k},e_{k}\right)=d\eta_{q}\left(E_{h}\left(q\right),E_{k}\left(q\right)\right)$,
for all $q\in M^{2n+1}$. Note that, by construction, we have
$\nabla_{E_h}E_k=0$ and $\nabla_{\xi}E_k=0$ for all
$h,k\in\left\{1,\ldots,2n\right\}$. From this and the expression
of the torsion of the bi-Legendrian connection (cf. $\S$
\ref{preliminari}), we get
\begin{gather}
    \left[Y_{i},Y_{j}\right]=\left[X_{i},X_{j}\right]=\left[Y_{i},\xi\right]=\left[X_{i},\xi\right]=0 \label{uno}\\
    \left[Y_{i},X_{j}\right]=-T\left(Y_i,X_{j}\right)=-2d\eta\left(Y_i,X_{j}\right){\xi}=\delta_{ij}\xi,\label{tre}
\end{gather}
for all $i,j\in\left\{1,\ldots,n\right\}$,  and
\eqref{uno}--\eqref{tre} imply that there exist local coordinates
$\{x_{1},\ldots,x_{n},$ $y_{1},\ldots,y_{n},z\}$ such that
$Y_{i}=\frac{\partial}{\partial y_{i}}$,
$X_{j}=\frac{\partial}{\partial x_{j}}+y_{j}\frac{\partial}{\partial
z}$, $\xi=\frac{\partial}{\partial z}$, for any
$i\in\left\{1,\ldots,n\right\}$. Note that from \eqref{relazioni2}
it follows that, with respect to these coordinates,
$d\eta=\sum_{i=1}^{n}dx_{i}\wedge dy_{i}$ from which we have
$d\left(\eta+\sum_{i=1}^{n}y_{i}dx_{i}\right)=0$ and so
$\eta=df-\sum_{i=1}^{n}y_{i}dx_{i}$, for some $f\in
C^{\infty}(M^{2n+1})$. But $\eta\left(Y_j\right)=0$,
$\eta\left(X_{j}\right)=0$ and $\eta\left(\xi\right)=1$ imply
$\frac{\partial f}{\partial y_j}=0$, $\frac{\partial f}{\partial
x_j}=0$ and $\frac{\partial f}{\partial z}=1$, respectively. So
$df=dz$ and, in this coordinate system we have that $T\cal F$ is
spanned by $Y_i=\frac{\partial}{\partial y_{i}}$, $T\cal G$ by
$X_i=\frac{\partial}{\partial
    x_{i}}+y_{i}\frac{\partial}{\partial z}$,
$i\in\left\{1,\ldots,n\right\}$, and
     the 1-form $\eta$ is given
    by $\eta=dz-\sum_{i=1}^{n}y_{i}dx_{i}$. Now we define a
tensor field $\phi$ and a Riemannian metric $g$ on $M^{2n+1}$
putting $\phi\xi=0$, $\phi Y_i =X_{i}$, $\phi X_{i}=-Y_i$, and
$g\left(Z,Z'\right)=-d\eta\left(Z,\phi Z'\right)$ for all
$Z,Z'\in\Gamma\left(\cal D\right)$,
$g\left(V,\xi\right)=\eta\left(V\right)$ for all
$V\in\Gamma(TM^{2n+1})$. A straightforward computation shows that
$\left(\phi,\xi,\eta,g\right)$ is indeed a Sasakian structure.
Finally, since, by construction,
$\nabla_{X_j}X_i=\nabla_{Y_j}X_i=\nabla_{\xi}X_i=0$,
$\nabla_{X_j}Y_i=\nabla_{Y_j}Y_i=\nabla_{\xi}Y_i=0$, we deduce
easily that $\nabla\phi=0$ and by Theorem \ref{tanakabilegendrian}
we get that $\nabla={^\ast\nabla}$.
\end{proof}

\begin{remark}
\emph{Assuming in Theorem \ref{sasaki1} and Theorem \ref{kahler1}
that the manifold is also simply connected we have that
$(M^{2n+1},\eta)$ and $(M^{2n},\omega)$ coincide  with
$\mathbb{R}^{2n+1}$ and $\mathbb{R}^{2n}$ with their usual contact
and symplectic structure, respectively.}
\end{remark}

Removing the initial hypothesis of $T\cal{G}$ being an affine
transversal distribution for $\cal F$, we have the following result.

\begin{theorem}\label{sasaki2}
Let $(M^{2n+1},\eta)$ be a contact manifold foliated by a flat
Legendrian foliation $\cal{F}$. Then there exists a Sasakian
structure $\left(\phi,\xi,\eta,g\right)$ on $(M^{2n+1},\eta)$ whose
Tanaka-Webster connection coincides with the bi-Legendrian
connection associated to the almost bi-Legendrian structure
$\left(L,Q\right)$, where $L=T\cal F$ and $Q=\phi L$.
\end{theorem}
\begin{proof}
In \cite{jayne0} it has been proved that given a flat Legendrian
foliation $\cal F$ of a contact manifold $(M^{2n+1},\eta)$, there
exists a canonical contact metric structure
$\left(\phi,\xi,\eta,g\right)$ such that
$(M^{2n+1},\phi,\xi,\eta,g)$ is a Sasakian manifold. This Sasakian
structure is defined in the following way. By the Darboux theorem
for Legendrian foliations (cf. \cite{pang}) for any point of
$M^{2n+1}$ there exists an open neighborhood with local coordinates
$\left\{x_1,\ldots,x_n,y_1,\ldots,y_n,z\right\}$ such that
$\eta=dz-\sum_{i=1}^{n}y_idx_i$, $\xi=\frac{\partial}{\partial z}$,
and $\cal F$ is locally spanned by the vector fields
$Y_i:=\frac{\partial}{\partial y_i}$,
$i\in\left\{1,\ldots,n\right\}$. Now consider the contact metric
structure $\left(\phi_U,\xi,\eta,g_U\right)$ on $U$ given by
\begin{equation*}
\phi_U=\left(%
\begin{array}{ccc}
  0 & I_n & 0 \\
  -I_n & 0 & 0 \\
  0 & Y & 0 \\
\end{array}%
\right), \ \
g_U=\left(%
\begin{array}{ccc}
  \delta_{ij}+y_iy_j & 0 & -y_i \\
  0 & \delta_{ij} & 0 \\
  -y_i & 0  & 1 \\
\end{array}%
\right),
\end{equation*}
where $Y$ is the $(1\times n)$-matrix given by
$Y=\left(y_1,\ldots,y_n\right)$. It is known (cf. \cite{yano1}) that
$\left(\phi_U,\xi,\eta,g_U\right)$ is a Sasakian structure on $U$.
Next, we consider an open covering of $M^{2n+1}$ by Darboux
neighborhoods as above, and using the fact that the leaves of
$\cal F$ have a natural flat affine structure it can be proved that these Sasakian structures fit together
for giving rise to a global Sasakian structure
$\left(\phi,\xi,\eta,g\right)$ on $M^{2n+1}$. Now consider the
conjugate Legendrian distribution $Q$ of $\cal F$, which by Lemma
\ref{flat} is also flat and which is generated by the vector
fields $X_i:=\phi Y_i=\frac{\partial}{\partial
x_i}+y_i\frac{\partial}{\partial z}$,
$i\in\left\{1,\ldots,n\right\}$. Applying \cite[Proposition
5.1]{cappellettimontano1}, we get
$\nabla_{X_j}X_i=\nabla_{Y_j}X_i=\nabla_{\xi}X_i=0$,
$\nabla_{X_j}Y_i=\nabla_{Y_j}Y_i=\nabla_{\xi}Y_i=0$, from which
$\nabla\phi=0$. Then, applying again Theorem
\ref{tanakabilegendrian}, we conclude that $\nabla$ coincides with
the Tanaka-Webster connection of
$(M^{2n+1},\phi,\xi,\eta,g)$.
\end{proof}

\begin{remark}
\emph{Note that the Legendrian distribution $Q$ of Theorem
\ref{sasaki2} is, a posteriori, integrable because of Theorem
\ref{tanakabilegendrian1}.}
\end{remark}

\begin{remark}
\emph{It should be noted that, by Corollary \ref{levicivita}, in
Theorem \ref{sasaki1} and \ref{sasaki2} the connections induced on
the leaves of $\cal F$ and $\cal G$ by the Levi Civita, the
Tanaka-Webster and the bi-Legendrian connection coincide.}
\end{remark}

\section{Examples and remarks}\label{esempi}

\begin{example}
\emph{Consider $\mathbb{R}^{2n+1}$ with its standard Sasakian
structure $\left(\phi,\xi,\eta,g\right)$ where}
\begin{equation*}
\eta=dz-\sum_{k=1}^{n}y_kdx_k \textrm{,  } \
\xi=\frac{\partial}{\partial z} \textrm{, } \
g=\eta\otimes\eta+\frac{1}{2}\sum_{k=1}^{n}\left(\left(dx_k\right)^2+\left(dy_k\right)^2\right)
\end{equation*}
\emph{and $\phi$ is represented by the $(2n+1)\times(2n+1)$
matrix}
\begin{equation*}
\left(%
\begin{array}{ccc}
  0 & I_n & 0 \\
  -I_n & 0 & 0 \\
  0 & y_1 \ \cdots \ y_n & 0 \\
\end{array}%
\right)
\end{equation*}
\emph{The standard bi-Legendrian structure $\left(L,Q\right)$ on \
$(\mathbb{R}^{2n+1},\phi,\xi,\eta,g)$ is given by $L$ \ $=$ \
$\textrm{span}\left\{X_1,\ldots,X_n\right\}$ and
$Q=\textrm{span}\left\{Y_1,\ldots,Y_n\right\}$, where, for all
$i\in\left\{1,\ldots,n\right\}$, $X_i:=\frac{\partial}{\partial
y_i}$ and $Y_i:=\frac{\partial}{\partial
x_i}+y_i\frac{\partial}{\partial z}$. It is easy to check that $\phi
X_i=Y_i$ for all $i\in\left\{1,\ldots,n\right\}$ and that $L$ and
$Q$ define two orthogonal flat Legendrian foliations on
$\mathbb{R}^{2n+1}$. Let $\nabla$ be the corresponding bi-Legendrian
connection. A straightforward computation shows that
$\nabla_{X_i}X_j=\nabla_{Y_i}X_j=\nabla_{\xi}X_j=0$ and
$\nabla_{X_i}Y_j=\nabla_{Y_i}Y_j=\nabla_{\xi}Y_j=0$. Using these
relations we have $\nabla\phi=0$ and so, by Proposition
\ref{metrica}, $\nabla g=0$. Then, by Theorem
\ref{tanakabilegendrian}, the bi-Legendrian connection $\nabla$
coincides with the Tanaka-Webster connection on
$(\mathbb{R}^{2n+1},\phi,\xi,\eta,g)$. In particular, with the
notation of Remark \ref{osservazione},
$L\in\Im_{\mathbb{R}^{2n+1}}$. Another consequence is that
 the  Tanaka-Webster connection on
$\mathbb{R}^{2n+1}$ is everywhere flat since $\nabla$ is flat (cf.
\cite{mino2})}.
\end{example}

\begin{corollary}\label{curvatura0}
Let $\cal{F'}$ be any Legendrian foliation on $\mathbb{R}^{2n+1}$
belonging to $\Im_{\mathbb{R}^{2n+1}}$. Then the curvature of the
corresponding bi-Legendrian connection vanishes identically.
\end{corollary}
\begin{proof}
$\cal{F'}$ is a flat Legendrian foliation on $\mathbb{R}^{2n+1}$
such that its conjugate Legendrian distribution is integrable
 and $\nabla'g=0$, where $\nabla'$ denotes the
bi-Legendrian connection associated to $\cal{F'}$. So, by Remark
\ref{osservazione} we have $\nabla=\nabla'$, $\nabla$ denoting the
bi-Legendrian connection associated to the standard bi-Legendrian
structure on $\mathbb{R}^{2n+1}$. In particular the curvature
tensor fields of the two connections must coincide and the result
follows from the flatness of $\nabla$.
\end{proof}
\\

Now we give an example of a Sasakian manifold endowed with a
non-flat bi-Legendrian structure for which the corresponding
bi-Legendrian connection is metric but does not coincide with the
Tanaka-Webster connection.

\begin{example}
\emph{Consider the sphere
$S^3=\left\{\left(x_1,x_2,x_3,x_4\right)\in\mathbb{R}^4:x_1^2+x_2^2+x_3^2+x_4^2=1\right\}$ with the following Sasakian structure:}
\begin{equation*}
\eta=x_3dx_1+x_4dx_2-x_1dx_3-x_2dx_4, \
\xi=x_3\frac{\partial}{\partial x_1}+x_4\frac{\partial}{\partial
x_2}-x_1\frac{\partial}{\partial x_3}-x_2\frac{\partial}{\partial
x_4},
\end{equation*}
\begin{equation*}
g=\left(%
\begin{array}{cccc}
  1 & 0 & 0 & 0 \\
  0 & 1 & 0 & 0 \\
  0 & 0 & 1 & 0 \\
  0 & 0 & 0 & 1 \\
\end{array}%
\right), \ \phi=\left(%
\begin{array}{cccc}
  0 & 0 & -1 & 0 \\
  0 & 0 & 0 & -1 \\
  1 & 0 & 0 & 0 \\
  0 & 1 & 0 & 0 \\
\end{array}%
\right).
\end{equation*}
\emph{Set $X:=x_2\frac{\partial}{\partial
x_1}-x_1\frac{\partial}{\partial x_2}-x_4\frac{\partial}{\partial
x_3}+x_3\frac{\partial}{\partial x_4}$ and $Y:=\phi
X=x_4\frac{\partial}{\partial x_1}-x_3\frac{\partial}{\partial
x_2}+x_2\frac{\partial}{\partial x_3}-x_1\frac{\partial}{\partial
x_4}$, and consider the $1$-dimensional distributions $L$ and $Q$ on
$S^3$ generated by $X$ and $Y$, respectively. An easy computation
shows that $\left[X,\xi\right]=-2Y$, $\left[Y,\xi\right]=2X$,
$\left[X,Y\right]=2\xi$. Thus $L$ and $Q$ defines two Legendrian
foliations on the Sasakian manifold $(S^3,\phi,\xi,\eta,g)$ which
are orthogonal and not flat. For the bi-Legendrian connection
corresponding to this bi-Legendrian structure, we have, after a
straightforward computation, $\nabla_X X=\nabla_X Y=\nabla_X
\xi=\nabla_Y X=\nabla_Y Y=\nabla_Y \xi=0$. Therefore $\nabla\phi=0$.
But $T\left(\xi,\phi V\right)=-\phi T\left(\xi,V\right)$ for all
$V\in\Gamma(TS^{3})$ is not satisfied; indeed $T\left(\xi,\phi
Y\right)=-T\left(\xi,X\right)=\left[\xi,X\right]=2Y$ and on the
other hand $\phi T\left(\xi,Y\right)=-\phi\left[\xi,Y\right]=2\phi
X=2Y$, so that $T\left(\xi,\phi Y\right)=-\phi T\left(\xi,X\right)$
holds if and only if $Y=0$.}
\end{example}

We conclude with an example of a bi-Legendrian structure on a
non-Sasakian manifold.

\begin{example}
\emph{Let \ $\mathfrak{g}$ \ be \ a \ $\left(2n+1\right)$-dimensional   \
 Lie \ algebra \ with basis \ $\{ X_{1},\ldots ,X_{n},$ \ $Y_{1},\ldots ,Y_{n},\xi\} $. The Lie bracket
is defined in the following way:}
\begin{gather*}
\left[ X_{i},X_{j}\right]=0  \text{ \emph{for any} }i,j\in \left\{
1,\ldots
,n\right\}, \ \left[ Y_{i},Y_{j}\right]=0  \text{ \emph{for any} }i\neq 2,\\
\left[ Y_{2},Y_{j}\right] =2Y_{j} \ \text{ \emph{for any} } j\neq
2, \ \left[X_{1},Y_{1}\right]=2{\xi }-2X_{2}, \
\left[X_{1},Y_{j}\right] =0  \text{ \emph{for any} }j\geq 2,\\
\left[ X_{h},Y_{k}\right]=\delta _{hk}\left( 2{\xi }%
-2X_{2}\right)  \text{ \emph{for any} }h,k\geq 3\text{,  } \
\left[X_{2},Y_{j}\right] =2X_{j}  \text{ \emph{for any} }j\neq 2\text{, }\\
\left[X_{2},Y_{2}\right]=2{\xi }\text{, } \
\left[ X_{k},Y_{1}\right]=\left[ X_{k},Y_{2}\right] =0 \ \text{ \emph{for any} }k\geq 3%
\text{,} \\
\left[ \xi,X_{j}\right] =0
\text{ \emph{and} } \ \left[ \xi,Y_{j}\right] =2X_{j} \ \text{
\emph{for any} }j\in \left\{ 1,\ldots ,n\right\} \text{,}
\end{gather*}
\emph{Let $G$ be a Lie group whose Lie algebra is $\mathfrak{g}$.
On $G$ one can define a contact metric structure by defining $\phi
\xi =0$, $\phi X_{i} =Y_{i}$, $\phi Y_{i} =-X_{i}$, for all $i\in
\left\{ 1,\ldots ,n\right\} $, considering the left invariant
Riemannian metric $g$ such that $\left\{ X_{1},\ldots
,X_{n},Y_{1},\ldots ,Y_{n},\xi\right\} $ is an orthonormal frame
and, finally, defining the 1-form $\eta$ as the dual 1-form of the
vector field $\xi$ with respect to the metric $g$. It can be
proved (cf. \cite{boeckx}) that $\left(G,\phi,\xi,\eta,g\right)$
is a contact $\left(\kappa,\mu\right)$-manifold with $\kappa=0$
and $\mu=4$ and so it is certainly non-Sasakian. Let $L$ and $Q$
be the $n$-dimensional distributions generated, respectively, by
$X_1,\ldots,X_n$ and $Y_1,\ldots,Y_n$. They can be viewed also as
the eigenspaces of the eigenvectors $\lambda$ and $-\lambda$ of
the operator $h$, where $\lambda=\sqrt{1-\kappa}=1$. As remarked
in Example \ref{kappamu}, $L$ and $Q$ define two orthogonal
Legendrian foliations of the contact metric manifold
$\left(G,\phi,\xi,\eta,g\right)$, and the corresponding
bi-Legendrian connection satisfies $\nabla g=0$, $\nabla\phi=0$.
Nevertheless it does not coincide with the Tanaka-Webster
connection of $\left(G,\phi,\xi,\eta,g\right)$. Indeed
$T\left(\xi,\phi
X_1\right)=-T\left(Y_1,\xi\right)=-\left[\xi,Y_1\right]_L=-2X_1$
and, on the other hand,
$T\left(\xi,X_1\right)=-T\left(X_1,\xi\right)=-\left[\xi,X_1\right]_Q=0$,
so $T\left(\xi,\phi X_1\right)\neq-\phi T\left(\xi,X_1\right)$.}
\end{example}

\section*{Appendix}
Recall that a Lagrangian foliation of a symplectic manifold
$(M^{2n},\omega)$ is an $n$-dimensional foliation $\cal
F$ of $M^{2n}$ such that $\omega\left(X,X'\right)=0$ for any
$X,X'\in\Gamma\left(T\cal F\right)$. A bi-Lagrangian structure on
$(M^{2n},\omega)$ is nothing but a pair of transversal
Lagrangian foliations $\left(\cal F,\cal G\right)$ on
$(M^{2n},\omega)$. In \cite{hess} H. Hess proved that, given two
transversal Lagrangian distributions $L$ and $Q$ on $M^{2n}$,
there exists a unique symplectic connection $\nabla$ on $M^{2n}$
preserving the distributions $L$ and $Q$ and whose torsion tensor
field satisfies
\begin{equation}\label{torsione}
T\left(X,Y\right)=0
\end{equation}
for all $X\in\Gamma\left(L\right)$ and $Y\in\Gamma\left(Q\right)$.
This connection is called the \emph{bi-Lagrangian connection}
associated to $\left(L,Q\right)$ and if $L$ and $Q$ are
integrable, i.e. if they define a bi-Lagrangian structure on
$M^{2n}$, $\nabla$ is torsion free and it is flat along the leaves
of the foliations. In this Appendix we prove the already stated
Theorem \ref{kahler1}, which, at the knowledge of the author, has
not been proved yet elsewhere.

\begin{lemma}[\cite{etayo1}]\label{metrica1}
Let $\left(\cal F,\cal G\right)$ be a bi-Lagrangian structure on
the symplectic manifold $(M^{2n},\omega)$. Let
$\left(J,\omega,g\right)$ be a Hermitian structure on
$(M^{2n},\omega)$. Then for the bi-Lagrangian
connection associated to $\left(\cal F,\cal G\right)$ we have
$\nabla g=0$ if and only if $\nabla J=0$.
\end{lemma}
\begin{proof}[Proof of Theorem \ref{kahler1}]
First of all note that, as in Theorem \ref{sasaki1}, the assumption
of $T\cal G$ being an affine transversal distribution implies that
$\nabla$ is everywhere flat. Fixed a point $x$ of $M^{2n}$, there
exists a basis
$\left\{e_{1},\ldots,e_{n},e_{n+1},\ldots,e_{2n}\right\}$ of
$T_{x}M^{2n}$ such that $\left\{e_{1},\ldots,e_{n}\right\}$ is a
basis of $T_{x}\cal F$, $\left\{e_{n+1},\ldots,e_{2n}\right\}$ is a
basis of $T_{x}\cal G$ and
\begin{equation}\label{relazionii}
\omega_x\left(e_{i},e_{j}\right)=\omega_x\left(e_{n+i},e_{n+j}\right)=0, \
\ \omega_x\left(e_{i},e_{n+j}\right)=-\frac{1}{2}\delta_{ij}
\end{equation}
for all $i,j\in\left\{1,\ldots,n\right\}$. For each
$k\in\left\{1,\ldots,2n\right\}$ we define a vector field $E_k$ on
$M^{2n}$ by the $\nabla$-parallel transport along curves of the
vector $e_k$. Note that, for all $y\in M^{2n}$,
$E_{k}\left(y\right)$ does not depend on the curve joining $x$ and
$y$, since $R\equiv 0$. Setting $X_i:=E_{n+i}$ and $Y_i:=E_i$, we
obtain $2n$ vector fields on $M^{2n}$ such that, for each
$i\in\left\{1,\ldots,n\right\}$, $Y_{i}\in\Gamma\left(T\cal
F\right)$ and $X_{i}\in\Gamma\left(T\cal G\right)$, because the
parallel transport preserves the distributions $T\cal{F}$ and
$T\cal{G}$. Moreover, since $\nabla\omega=0$, \eqref{relazionii}
hold at any point of $M^{2n}$, that is
\begin{equation}\label{relazionii2}
\omega_y\left(Y_{i}\left(y\right),Y_{j}\left(y\right)\right)=\omega_y\left(X_{i}\left(y\right),X_{j}\left(y\right)\right)=0,
\     \
\omega_y\left(Y_{i}\left(y\right),X_{j}\left(y\right)\right)=-\frac{1}{2}\delta_{ij}
\end{equation}
for any $y\in M^{2n}$ and $i,j\in\left\{1,\ldots,n\right\}$. Note
that, by construction, we have $\nabla_{E_h}E_k=0$ for all
$h,k\in\left\{1,\ldots,2n\right\}$. From this and \eqref{torsione}
we get
\begin{equation}\label{uno1}
    \left[Y_{i},Y_{j}\right]=\left[X_{i},X_{j}\right]=\left[Y_{i},X_{j}\right]=0
\end{equation}
for all $i,j\in\left\{1,\ldots,n\right\}$,  and \eqref{uno1} imply
the existence of coordinates
$\{x_{1},\ldots,x_{n},y_{1},\ldots,y_{n}\}$ such that for each
$i\in\left\{1,\ldots,n\right\}$ $Y_{i}=\frac{\partial}{\partial
y_{i}}$ and $X_{j}=\frac{\partial}{\partial x_{j}}$. So in this
coordinate system we have that $T\cal F$ is spanned by
$Y_i=\frac{\partial}{\partial y_{i}}$, $T\cal G$ by
$X_i=\frac{\partial}{\partial
    x_{i}}$,
$i\in\left\{1,\ldots,n\right\}$, and, moreover, by
\eqref{relazionii2}, $\omega=\sum_{i=1}^{n}dx_{i}\wedge dy_{i}$. Now
we define a tensor field $J$ and a Riemannian metric $g$ on $M^{2n}$
putting, for each $i\in\left\{1\,\ldots,n\right\}$, $J Y_i =X_{i}$,
$J X_{i}=-Y_i$, and $g\left(V,W\right)=-\omega\left(V,J W\right)$
for all $V,W\in\Gamma(TM^{2n})$. A straightforward computation shows
that $\left(J,\omega,g\right)$ is indeed a K\"{a}hlerian structure.
Finally, since, by construction,
$\nabla_{X_j}X_i=\nabla_{Y_j}X_i=\nabla_{X_j}Y_i=\nabla_{Y_j}Y_i=0$,
we deduce easily that $\nabla J=0$, which, by Lemma \ref{metrica1},
imply $\nabla g=0$. Thus $\nabla$ coincides with the Levi Civita
connection of $(M^{2n},J,\omega,g)$.
\end{proof}

\end{document}